\documentclass{article}
\usepackage[english]{babel}
\usepackage[utf8]{inputenc}

\usepackage{bm}
\usepackage{amsmath}
\usepackage{graphicx}
\usepackage{color}
\usepackage{multirow}
\usepackage{subcaption}
\usepackage{amsfonts}
\newcommand{\const}{\mathop{\rm const}\nolimits}

\graphicspath{ {./figs/} } % папки с рисунками

\usepackage{placeins}
\usepackage{booktabs, multirow} % for borders and merged ranges
\usepackage{soul}% for underlines
\usepackage[table]{xcolor} % for cell colors
\usepackage{changepage,threeparttable} % for wide tables
\date{}

\title{Numerical investigation and factor analysis of the spatial-temporal multi-species competition problem}

\author{
Maria Vasilyeva 
\thanks{Department of Mathematics and Statistics, Texas A\&M University  - Corpus Christi,   Corpus Christi, Texas, USA. Email: {\tt maria.vasilyeva@tamucc.edu}.}
\and
Youwen Wang
\thanks{Department of Mathematics and Statistics, Texas A\&M University  - Corpus Christi,   Corpus Christi, Texas, USA. Email: {\tt ywang23@islander.tamucc.edu}.}
\and
Sergei Stepanov
\thanks{Laboratory of Computational Technologies for Modeling Multiphysical and Multiscale Permafrost Processes, North-Eastern Federal University, Yakutsk, Republic of Sakha (Yakutia), 677980, Russia}
\and
Alexey Sadovski
\thanks{Department of Mathematics and Statistics, Texas A\&M University  - Corpus Christi,   Corpus Christi, Texas, USA. Email: {\tt alexey.sadovski@tamucc.edu}.}
}

\begin{document}

\maketitle

\begin{abstract}
In this work, we consider the spatial-temporal multi-species competition model.  A mathematical model is described by a coupled system of nonlinear diffusion-reaction equations.  We use a finite volume approximation with semi-implicit time approximation for the numerical solution of the model with corresponding boundary and initial conditions.  
To understand the effect of the diffusion to solution in one and two-dimensional formulations,  we present numerical results for several cases of the parameters related to the survival scenarios.  The random initial conditions' effect on the time to reach equilibrium is investigated.  The influence of diffusion on the survival scenarios is presented. 
In real-world problems, values of the parameters are usually unknown and vary in some range.  In order to evaluate the impact of parameters on the system stability, we simulate a spatial-temporal model with random parameters and perform factor analysis for two and three-species competition models.  
\end{abstract}

\section{Introduction}

%APPLICATION
Understanding the stability of ecosystems is of fundamental importance to ecology \cite{okuboDiffusionEcologicalProblems2001, murray2002mathematical, marchukmathematical}.  
Fundamental mathematical models of such systems are described by a coupled system of ordinary differential equations (ODEs). 
The Lotka–Volterra Competition model (LVC)  is a  basic model that describes the dynamics of the species population competing for some shared resource. 
The LVC model has been applied in many areas, including biological systems, industry and economics \cite{marasco2016market, zhang2013maritime, windarto2020modification, wang2021competition, khan2019fractional}.   
For example,  the model can be used to simulate the marsh ecosystems for the wetlands at the Nueces River mouth 
\cite{montagnaModelingEffectWater2017}.

The LVC model is based on the logistic population model
\[
\frac{d u}{d t} = r u (1 - u),
\]
where $u$ is the size of the population at a given time $t$and $r > 0$ is per-capita growth rate.  

For the general case of multispecies competition, we have
\[
\begin{split}
\frac{d u^{(k)}}{d t} = 
 r_k   u^{(k)} (1 - u^{(k)})
 - \sum_{l \neq k} \alpha_{kl} u^{(k)} u^{(l)}, \quad \forall k = 1,..,M,
\end{split}
\]
where 
$u^{(k)}$ is the population of the $k$-th species, 
$r_k$ is the $k$-th population growth rate,  
$\alpha_{kl}$ is the  interaction coefficient due to competition ($u^{(l)}$ compete with $u^{(k)}$),  and 
$M$ is the number of species (equations).

% ----- diffusion -----
In such systems, a temporal dynamic model can represent and describe the behaviors of the entire system.   However, real ecosystems interact in different locations, and spatial structure is important and has a large impact on the final equilibrium state \cite{chenSpatiotemporalEcologicalModels2011, zelnikImpactSpatialTemporal2018}.
Systems of the partial differential equations (PDEs) are used to describe spatial-temporal systems.  A mathematical model is described by a coupled system of unsteady nonlinear reaction-diffusion equations. 
For the  multispecies interaction, we have 
\[
\begin{split}
\frac{\partial u^{(k)}}{\partial t} 
- \nabla \cdot (\varepsilon_k \nabla u^{(k)}) = 
 r_k   u^{(k)} (1 - u^{(k)})
 - \sum_{l \neq k} \alpha_{kl} u^{(k)} u^{(l)},
\end{split}
\]
where  
$\varepsilon_k$ is the diffusion coefficient. 
This model incorporates spatial structure by adding diffusion terms in equations and considering the system of equations in domain $\Omega \subset R^d$, where $d$ is the spatial dimension.

% METHOD 
% ----- cite paper for making parameters as function -----
Previous research has put much effort into modifying the Lotka–Volterra competition model (LVC). The parameters $r_k$, $\alpha_{kl}$, and $\varepsilon_k$ can be modeled as functions instead of constant. Diffusion $\varepsilon_k$, more specifically, traveling waves and their minimal-speed selection mechanisms in the LV model can be studied by applying the upper-lower solution technique on the cooperative system \cite{alhasanatMinimalspeedSelectionTraveling2019}. The crowding effect of diffusion $\varepsilon_k$ can also be modeled within LVC model \cite{gavinaMultispeciesCoexistenceLotkaVolterra2018}. Advection rate can be introduced to the LVC model as a supplement of diffusion to bring much richer phenomena \cite{zhouLotkaVolterraCompetitionSystem2016,zhaoLotkaVolterraCompetition2016}. Carrying capacity can be replaced by function instead of constant and can vary between species; connection topology can be modeled as ‘loop’, ‘star’, ‘chain’, or ‘full’ connection when there are more than two species in the system \cite{dakosIdentifyingBestindicatorSpecies2018}. 
In addition to making parameters endogenous, some other efforts have also been made to modify the LVC model to simulate real-world problems. For example, simulation adding small immigration into the prey or predator population can stabilize the LV system \cite{taharaAsymptoticStabilityModified2018}. The random fluctuating environment can also be modeled within the LV model to show how switching between environments can make survival harder \cite{benaimLotkaVolterraRandomly2016}. Stochastic noises can be introduced to the LV population model and are presented to play an essential role in the permanence and characterization of the system \cite{liuPermanenceStochasticLotka2017}.

% ----- cite paper for using different methods of solving LVC model -----
Multiple methods are used for solving the LVC model, such as Haar wavelet (HW), Adams-Bashforth-Moulton (ABM) methods \cite{kumarStudyFractionalLotkaVolterra2020}, and finite element method. 
% ----- cite paper for analyzing result of simulated data -----
After simulation, multiple methods have been used to analyze the simulated species interaction data. Real-world data can be used to evaluate the multi-species interaction model \cite{devarajanMultiSpeciesOccupancy2020}. For example, the population dynamics model of individual reefs can be compared with data of coral reefs in Pilbara \cite{boschettiSettingPrioritiesConservation2020}. Numerical models and machine learning can be combined to identify the factors that influence Alexandrium catenella blooms \cite{baekIdentificationInfluencingFactors2021}.

%CURRENT WORK
% ----- in this work -----
In this work, we consider a  spatial-temporal multi-species competition model in one- and two - dimensional formulations. For the numerical solution of the model with corresponding boundary and initial conditions, we construct a finite volume approximation with a semi-implicit time approximation.   
When two or more species compete for the same limited food source or in some way influence each other's growth, one or several of the species usually becomes extinct. 
To understand the diffusion effect of the solution,  we perform a numerical investigation for several cases of the parameters related to the survival scenarios. The influence of the diffusion and initial conditions on the survival scenarios is presented for two and three-species competition models. After understanding the effect of parameters and initial conditions on the test problems, we evaluate the impact of parameters on the system stability by simulating the spatial-temporal model with random input parameters and performing factor analysis.    
Simulations with different combinations of parameters and initial conditions support the hypothesis that such systems reach equilibrium sooner or later. 
It turns out that equilibrium depends only on system parameters (birth rate, competition, predation, etc.) and does not depend on the initial conditions. The time it takes for the system to be in the steady state depends on how far initial populations of species are located from the equilibrium point in the phase space of the model.

%STRUCTURE OF THIS PAPER
% in this work
The paper is organized as follows.  
Section 2 describes the mathematical model with fine-scale approximation using the finite volume method and a semi-implicit scheme for time approximation.  
In Section 3, we present numerical results for two and three-species competition models in one and two-dimensional formulations.  
In Section 4, we simulate a spatial-temporal model with random parameters and perform factor analysis to evaluate the impact of parameters on the system stability.  
The paper ends with a conclusion.

\section{Mathematical model with approximation by space and time}

The mathematical model is described by a coupled system of nonlinear diffusion - reaction equations in domain $\Omega \subset R^d$ ($d=1,2$) 
\begin{equation}
\label{m}
\begin{split}
\frac{\partial u^{(k)}}{\partial t} 
- \nabla \cdot (\varepsilon_k \nabla u^{(k)}) = 
 r_k   u^{(k)} (1 - u^{(k)})
 - \sum_{l \neq k} \alpha_{kl} u^{(k)} u^{(l)}, 
\quad x \in \Omega,  \quad t > 0,
\end{split}
\end{equation}
where 
$k = 1,..,M$, where $M$ is the number of species (equations). 
Here 
$u^{(k)}$ is the population of the $k$-th species, 
$r_k$ is the $k$-th population reproductive growth rate, 
$\alpha_{kl}$ is the  interaction coefficient due to competition ($u^{(l)}$ compete with $u^{(k)}$) and $\varepsilon_k$ is the diffusion coefficient c\cite{chairez2020spatial}. 

The system of equation is considered with initial conditions
\begin{equation}
\label{m-ic}
u^{(k)} = u^{(k)}_0, \quad x \in \Omega, \quad t = 0,
\end{equation}
and boundary conditions
\begin{equation}
\label{m-bc}
u^{(k)} = 0,  \quad x \in \partial \Omega,  \quad t > 0. 
\end{equation}

In this work, we consider following special cases:
\begin{itemize}
\item Two-species competition
\[
\begin{split}
\frac{\partial u^{(1)}}{\partial t} 
- \nabla \cdot (\varepsilon_1 \nabla u^{(1)}) & 
=  r_1  u^{(1)} (1 - u^{(1)}) 
 - \alpha_{12} u^{(1)} u^{(2)}, 
\\
\frac{\partial u^{(2)}}{\partial t} 
- \nabla \cdot (\varepsilon_2 \nabla u^{(2)}) & 
=  r_2  u^{(2)} (1 - u^{(2)}) 
 - \alpha_{21} u^{(1)} u^{(2)}, 
\end{split}
\]
\item Three-species competition
\[
\begin{split}
\frac{\partial u^{(1)}}{\partial t} 
- \nabla \cdot (\varepsilon_1 \nabla u^{(1)}) & 
=  r_1  u^{(1)} (1 - u^{(1)}) 
 - \alpha_{12} u^{(1)} u^{(2)}
 - \alpha_{13} u^{(1)} u^{(3)}, 
\\
\frac{\partial u^{(2)}}{\partial t} 
- \nabla \cdot (\varepsilon_2 \nabla u^{(2)}) & 
=  r_2  u^{(2)} (1 - u^{(2)}) 
 - \alpha_{21} u^{(1)} u^{(2)}
 - \alpha_{23} u^{(3)} u^{(2)}, 
\\
\frac{\partial u^{(3)}}{\partial t} 
- \nabla \cdot (\varepsilon_3 \nabla u^{(3)}) & 
=  r_3  u^{(3)} (1 - u^{(3)}) 
 - \alpha_{31} u^{(1)} u^{(3)} 
 - \alpha_{32} u^{(2)} u^{(3)}, 
\end{split}
\]
\end{itemize}

%%% FVM
For numerical solution of the problem \eqref{m} with initial and boundary conditions \eqref{m-ic}-\eqref{m-bc}, we construct the structured grid for domain $\Omega = [0,L]^d$ ($d = 1,2$) with $h = L/(N-1)$, where  $N$ is the number of nodes in each direction \cite{samarskii2001theory}.    
We set
\[
\frac{1}{|K_i|}\int_{K_i} u^{(k)} dx = u_i^{(k)},
\]
where $|K_i|$ is cell volume and $u_i^{(k)}$ is the cell average value of the function $u^{(k)}$ on cell $K_i$.

For the finite volume method, we have  
\begin{equation}
\frac{ \partial u_i^{(k)} }{\partial t} |K_i| 
+  \sum_{j} T_{k,ij} ( u^{(k)}_i - u^{(k)}_j ) 
= 
r_{k}   u_i^{(k)} (1 - u_i^{(k)}) |K_i|
 - \sum_{l \neq k} \alpha_{kl} u_i^{(k)}  u_i^{(l)} |K_i|, 
\end{equation}
with 
\[
 T_{k,ij} = \varepsilon_{k} \ |e_{ij}| / d_{ij}, 
\]
where $d_{ij}$ is the distance between to cell center points $x_i$ and $x_j$,  
$|e_{ij}|$ is the length of the interface between two cells $K_i$ and $K_j$. 
Note that, for structured uniform grid we have $d_{ij} = h$,  $|e_{ij}| = 1$  for one-dimensional case and $|e_{ij}|=h$ for two-dimensional problem. 

For the time approximation, we use a semi-implicit scheme \cite{samarskii1995computational, vabishchevich2013additive,afanasyeva2012unconditionally}
\begin{equation}
\frac{ u^{(k)}_i - \check{u}^{(k)}_i}{\tau} |K_i| +
\sum_{j} T_{k,ij} ( u^{(k)}_i - u^{(k)}_j ) = 
r_{k}   \check{u}_i^{(k)} (1 - \check{u}_i^{(k)}) |K_i|
 - \sum_{l \neq k} \alpha_{kl} \check{u}_i^{(k)}  \check{u}_i^{(l)} |K_i|, 
\end{equation}
where
$\tau$ is the time step size and $\check{u}^k_i$ is the solution from previous time layer. 
Note that, we obtain uncoupled system of equations and can solve equation for each component separately.

\section{Numerical results}

Before we make factor analysis for randomly generated parameters, we consider some special cases with parameters (growth rate "r", competition efficiency "$\alpha$", and diffusion rate "$\varepsilon$" ) fixed as constant and examine results such that either one species survives, two species survive, or three species survive.  

Let 
\[
r = (r_1, ..., r_M),  \quad 
\varepsilon = (\varepsilon_1, ..., \varepsilon_M),   \quad 
\alpha = \begin{pmatrix}
\alpha_{11} & \alpha_{12} & ... & \alpha_{1M} \\ 
\alpha_{21} & \alpha_{22} & ... & \alpha_{2M} \\ 
... & ... & ... & ... \\ 
\alpha_{M1} & \alpha_{M2} & ... & \alpha_{MM} 
\end{pmatrix}
\]

We present numerical results for multispecies competition in domain $\Omega$. We consider the following domain and boundary conditions configurations:
\begin{small}
\begin{itemize}
\item \textit{1D}: $\Omega = [0,1]$ with zero (fixed) boundary conditions  on $\partial \Omega$. 
\item \textit{2D(a)}: $\Omega = [0,1]^2$ with zero (fixed) boundary conditions on $\partial \Omega$. 
\item \textit{2D(b)}: $\Omega = [0,1]^2$ with zero (fixed) boundary conditions on left and right boundaries, and zero flux (free) boundary conditions on top and bottom boundaries
\end{itemize}
\end{small}
%The \textit{2D(a)} spatial boundary condition can be seen as an approximation of lake, while the \textit{2D(b)} spatial boundary condition can be seen as an approximation of river. 
In simulations, we use grid with $100$ nodes for one-dimensional problem and $25 \times 25$ grid for two-dimensional case. 
We simulate with $\tau = 1$ and set initial conditions $u_0^{(1)} =u_0^{(2)}= 0.5$ for two-species and $u_0^{(1)} =u_0^{(2)}=u_0^{(3)} = 0.5$ for three-species competition models.   
As of diffusion rate, we consider cases with regular diffusion $\varepsilon = D$ (which is between 0.01 and 0.1, and has same scale as other parameters) and small  diffusion $\varepsilon = D/10$. 

To represent result and compare final equilibrium state, we calculate average solution over computational domain $\Omega$ for each species
\[
\bar{u}^{(k)}(t) = \frac{1}{|\Omega|} \int_{\Omega} u^{(k)}(x, t) \ dx, 
\]
where $|\Omega|$ is the volume of domain $\Omega$. 
In most results, we perform simulations till both population reach equilibrium, $|\bar{u}^{(k)} - \check{\bar{u}}^{(k)}| < \epsilon$ with $\epsilon = 10^{-5}$ for each $k$.  

First,  we present results for the two-species interaction problem, where we simulate two sets of parameters related to the species survival: Case 1 (one species survive) and Case 2 (both species survive). We compare the dynamics of the average solution and final state for three configurations: 1D, 2D(a), and 2D(b). Furthermore, the investigation of the diffusion parameters scale to the final state is presented.
Next, we consider the three-species competition model with three cases of the test parameters: Case 1 (one species survive), Case 2 (two species survive), and Case 3 (all species survive). Similarly to the previous problems, we investigate the final state and dynamic for 1D, 2D(a), and 2D (b) for small and regular diffusion. 

After that, we present results for random diffusion to investigate the effect of the species diffusion coefficient on the final equilibrium state, survival group, and time to reach an equilibrium state.   Next, we consider the influence of the initial conditions on the equilibrium state and the time to reach it.  
%Finally, we simulate a two- and three-species competition model with randomly generated parameters and perform factor analysis. 

\subsection{Two-species competition}

We consider  two-species competition model and simulate two cases of the parameters:
\begin{itemize}
\item Case 1 (one species survive)
\[
D = (0.035, 0.014), \quad 
r = (0.074, 0.084), \quad
\alpha =
\begin{pmatrix}
0.0 & 0.074\\ 
0.013 & 0.0
\end{pmatrix}
\]
\item Case 2 (both species survive)
\[
D = (0.016, 0.014), \quad 
r = (0.083, 0.081), \quad
\alpha =
\begin{pmatrix}
0.0 & 0.053\\ 
0.049 & 0.0
\end{pmatrix}
\]
\end{itemize}
As of diffusion rate, we set $\varepsilon = D$ and ten times smaller  diffusion $\varepsilon = D/10$.

%\FloatBarrier
\begin{figure}[h!]
\centering
\begin{subfigure}{0.32\textwidth}
\centering
1D\\
\includegraphics[width=1\linewidth]{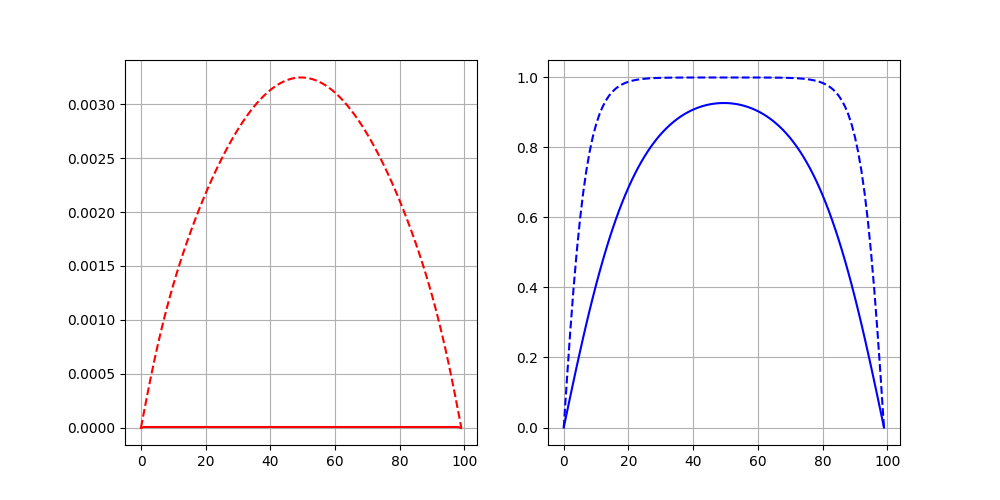}
\end{subfigure}
\begin{subfigure}{0.32\textwidth}
\centering
2D (a)\\
\includegraphics[width=1\linewidth]{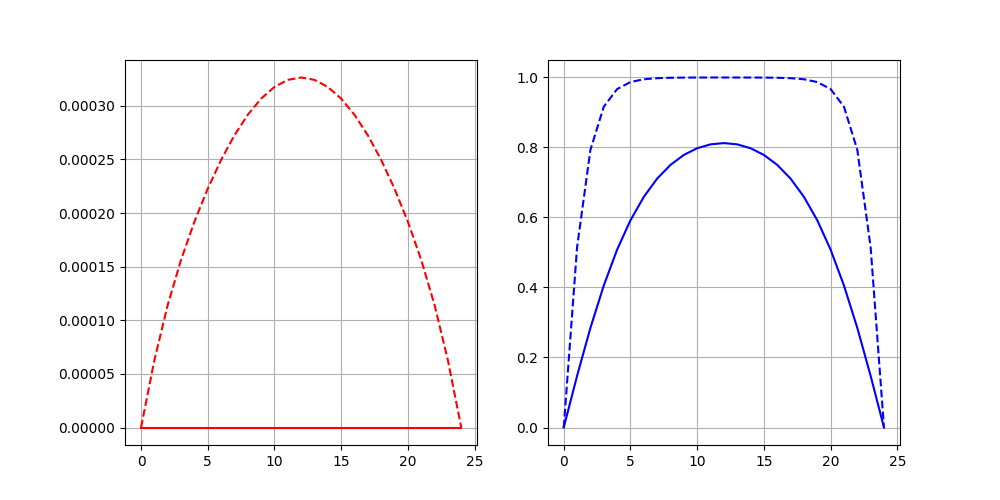}
\end{subfigure}
\begin{subfigure}{0.32\textwidth}
\centering
2D (b)\\
\includegraphics[width=1\linewidth]{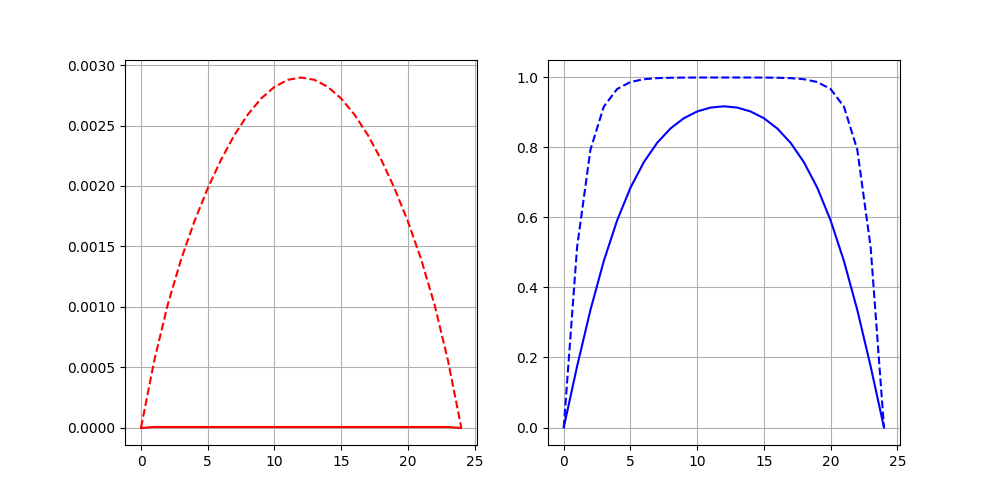}
\end{subfigure}\\
Case 1 (one species survive)
\\
\begin{subfigure}{0.32\textwidth}
\centering
\includegraphics[width=1\linewidth]{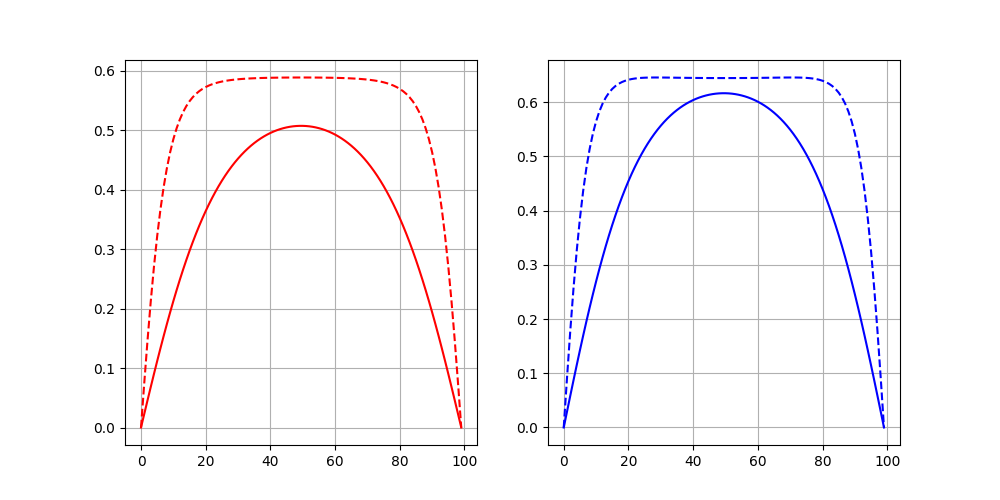}
\end{subfigure}
\begin{subfigure}{0.32\textwidth}
\centering
\includegraphics[width=1\linewidth]{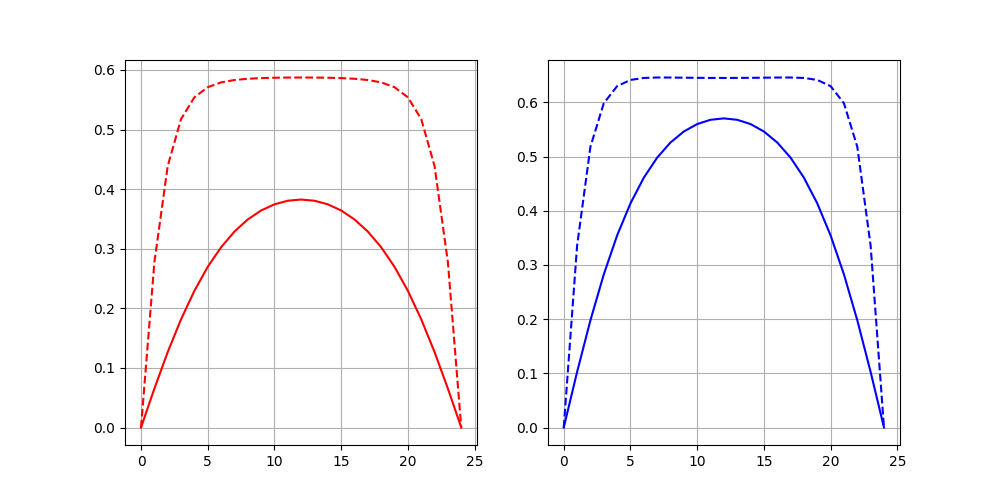}
\end{subfigure}
\begin{subfigure}{0.32\textwidth}
\centering
\includegraphics[width=1\linewidth]{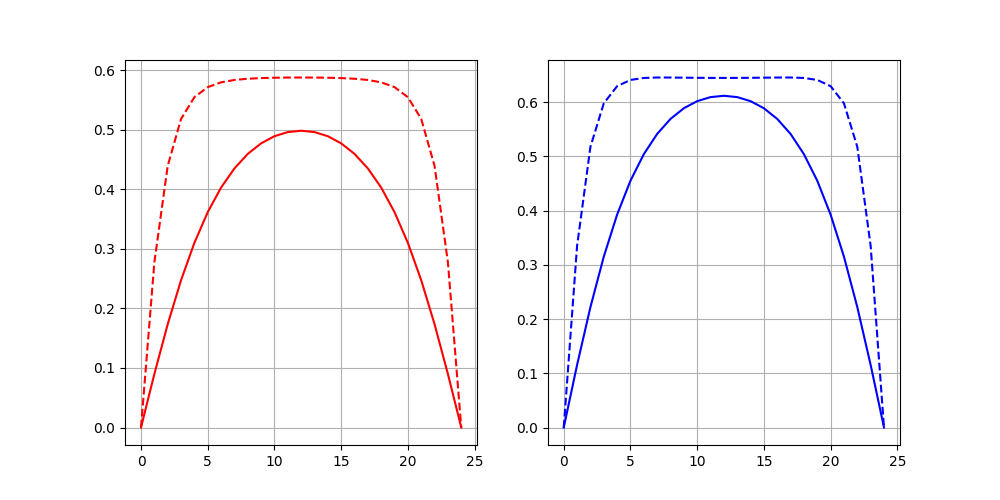}
\end{subfigure}\\
Case 2 (two species survive)
\caption{Solution at final time for regular diffusion $\varepsilon =  D$ (solid line) and small diffusion $\varepsilon = D/10$ (dashed line). 
Red color: first species. Blue color: second species }
\label{cc-1}
\end{figure}
%\FloatBarrier

\begin{figure}[h!]
\centering
\begin{subfigure}{0.32\textwidth}
\centering
1D\\
\includegraphics[width=0.8\linewidth]{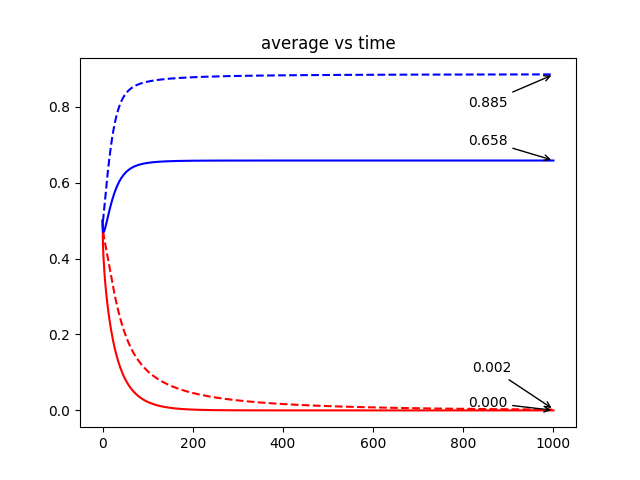}
\end{subfigure}
\begin{subfigure}{0.32\textwidth}
\centering
2D (a)\\
\includegraphics[width=0.8\linewidth]{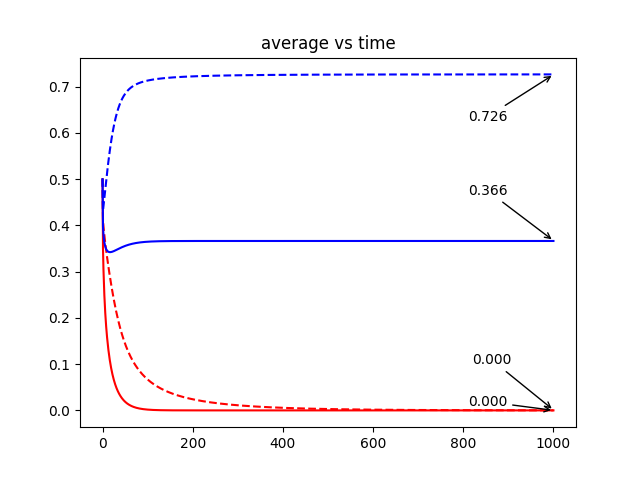}
\end{subfigure}
\begin{subfigure}{0.32\textwidth}
\centering
2D (b)\\
\includegraphics[width=0.8\linewidth]{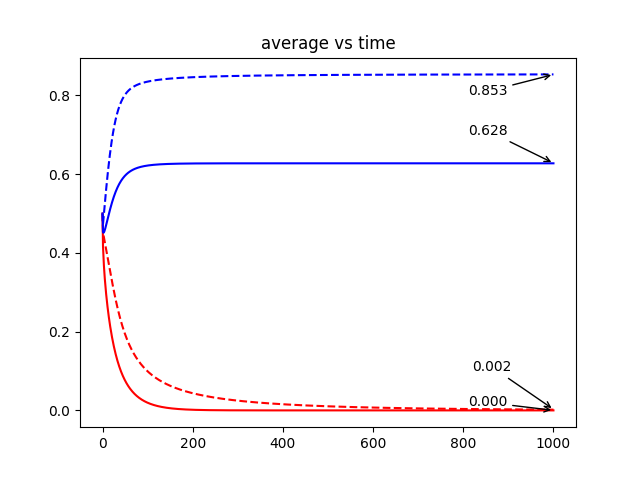}
\end{subfigure}\\
Case 1 (one species survive)
\\
\begin{subfigure}{0.32\textwidth}
\centering
\includegraphics[width=0.8\linewidth]{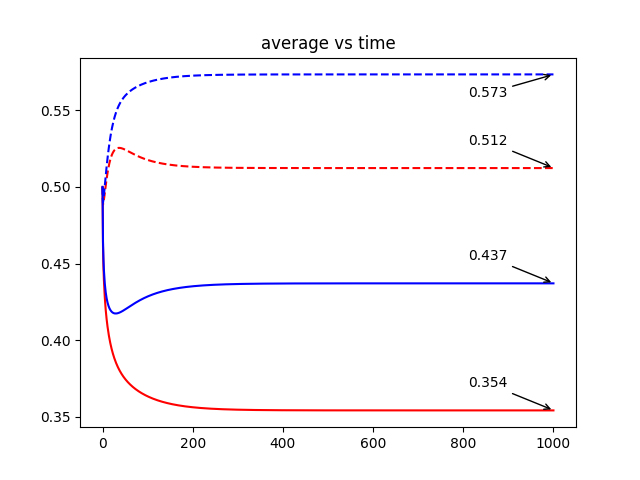}
\end{subfigure}
\begin{subfigure}{0.32\textwidth}
\centering
\includegraphics[width=0.8\linewidth]{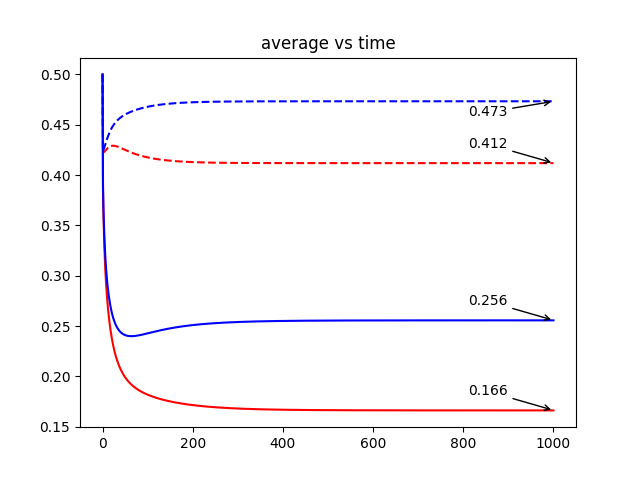}
\end{subfigure}
\begin{subfigure}{0.32\textwidth}
\centering
\includegraphics[width=0.8\linewidth]{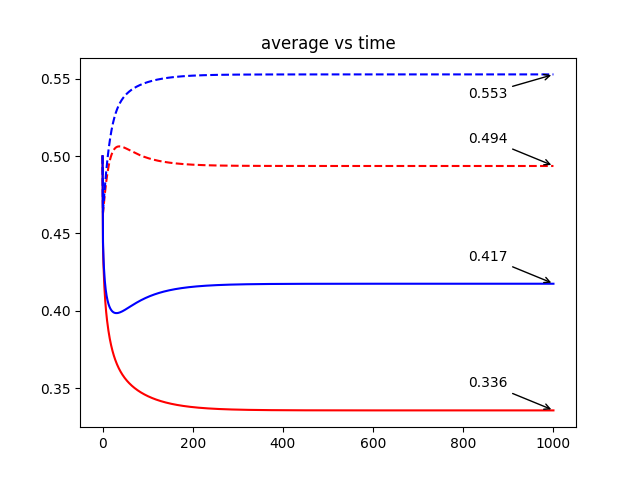}
\end{subfigure}\\
Case 2 (two species survive)
\caption{Dynamic of the average solution over domain for regular diffusion $\varepsilon =  D$ (solid line) and small diffusion $\varepsilon = D/10$ (dashed line). 
Red color: first species. Blue color: second species }
\label{cc-1t}
\end{figure}
%\FloatBarrier

%\FloatBarrier
\begin{figure}[h!]
\centering
\begin{subfigure}{0.49\textwidth}
\centering
2D (a)
\end{subfigure}
\begin{subfigure}{0.49\textwidth}
\centering
2D (b)
\end{subfigure}
\\
\begin{subfigure}{0.49\textwidth}
\centering
\includegraphics[width=1\linewidth]{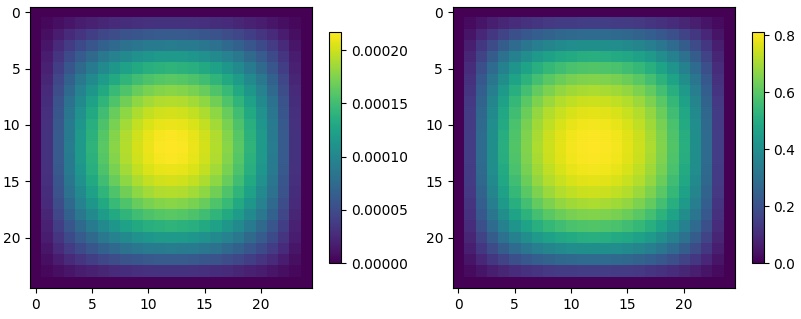}
\end{subfigure}
\begin{subfigure}{0.49\textwidth}
\centering
\includegraphics[width=1\linewidth]{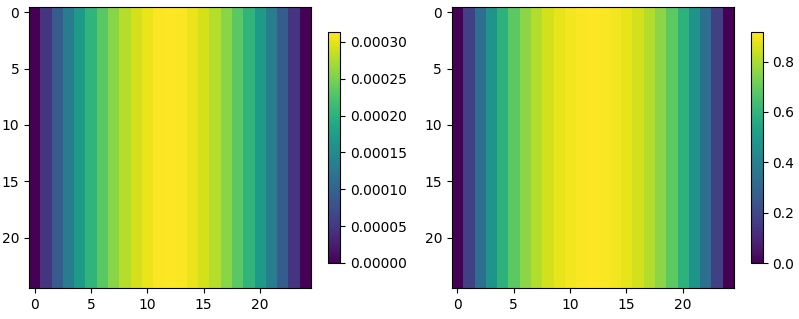}
\end{subfigure}\\
Case 1 (one species survive)
\\
\begin{subfigure}{0.49\textwidth}
\centering
\includegraphics[width=1\linewidth]{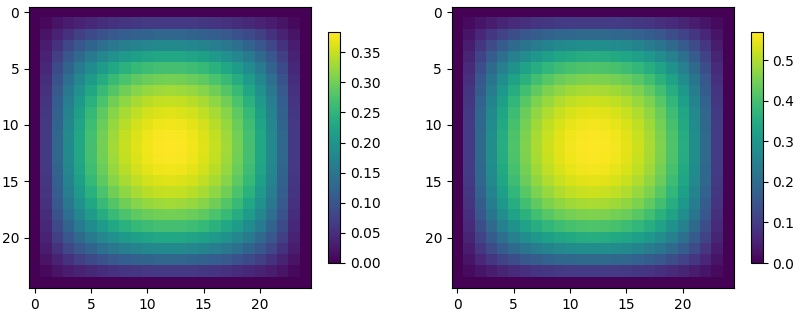}
\end{subfigure}
\begin{subfigure}{0.49\textwidth}
\centering
\includegraphics[width=1\linewidth]{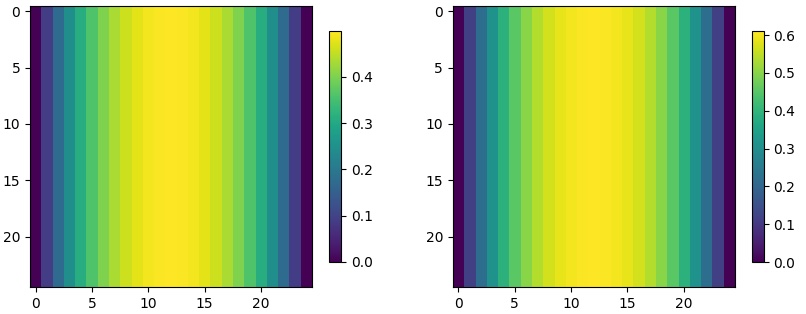}
\end{subfigure}\\
Case 2 (two species survive)
\caption{Solution at final time for regular diffusion $\varepsilon =  D$. 
First picture: first species. Second picture: second species }
\label{cc-2d}
\end{figure}
%\FloatBarrier

%\FloatBarrier
\begin{figure}[h!]
\centering
\begin{subfigure}{0.49\textwidth}
\centering
2D (a)
\end{subfigure}
\begin{subfigure}{0.49\textwidth}
\centering
2D (b)
\end{subfigure}
\\
\begin{subfigure}{0.49\textwidth}
\centering
\includegraphics[width=1\linewidth]{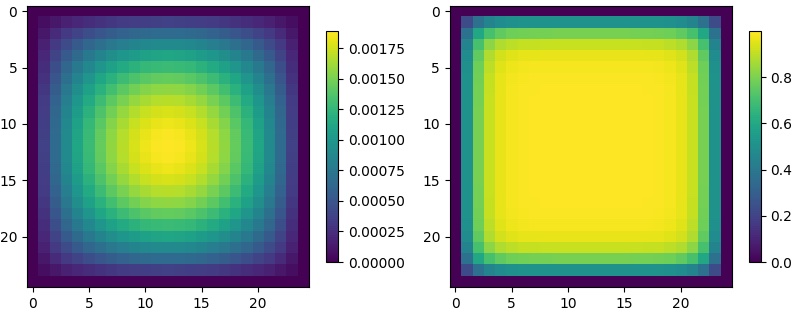}
\end{subfigure}
\begin{subfigure}{0.49\textwidth}
\centering
\includegraphics[width=1\linewidth]{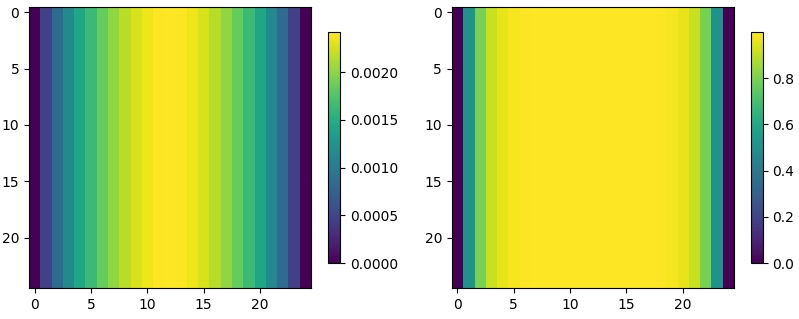}
\end{subfigure}\\
Case 1 (one species survive)
\\
\begin{subfigure}{0.49\textwidth}
\centering
\includegraphics[width=1\linewidth]{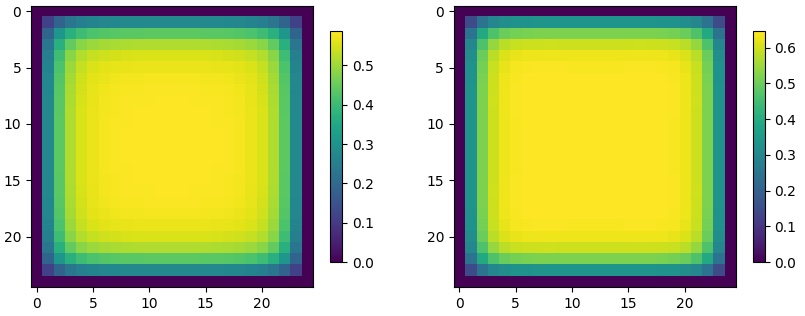}
\end{subfigure}
\begin{subfigure}{0.49\textwidth}
\centering
\includegraphics[width=1\linewidth]{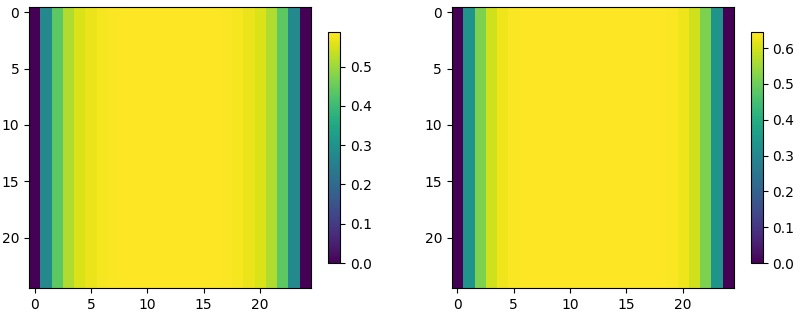}
\end{subfigure}\\
Case 2 (two species survive)
\caption{Solution at final time for  small diffusion $\varepsilon = D/10$. 
First picture: first species. Second picture: second species }
\label{cc-2ds}
\end{figure}
%\FloatBarrier

In Figures \ref{cc-1}, we plot the solution for one and two-dimensional formulations at the final time. Note that for 2d we plot solution over middle line ($y = 0.5$). We observed that the effect of boundary constraint (set at zero) is more severe in regular diffusion groups, causing a lower final population density than in small diffusion groups. With the same survival status (one species survives, or both species survive) and  boundary conditions (\textit{1D}, \textit{2D(a)}, or \textit{2D(b)}), the regular diffusion group always arrive at a lower population density equilibrium, compare to the small diffusion group. In other words, when the diffusion is smaller, the species can reach a higher population density equilibrium.

In Figures \ref{cc-1t}, we present the dynamic of the solution average over time for two species systems. We observed that the time to reach equilibrium is different among different diffusion conditions (regular diffusion $D$, or small diffusion $D/10$), survival status (one species survives, or both species survive) and boundary conditions (\textit{1D}, \textit{2D(a)}, or \textit{2D(b)}). In general, \textit{1D} and \textit{2D(b)} give similar solutions, while \textit{2D(a)} gives different solution. %This suggests that we can use \textit{1D} simulation to approximate \textit{2D(b)} condition, but we cannot use the former to approximate \textit{2D(a)} condition.

In the case where only one species survives, under \textit{1D} or \textit{2D(b)} boundary condition, for the species that survive, it takes the small diffusion group less time to reach a higher final population density equilibrium. In comparison, it took the regular diffusion group more time to reach a lower population density equilibrium above the initial population density. In contrast, under \textit{2D(a)}  boundary condition, for the species that survive, it takes the small diffusion group more time to reach a higher population density equilibrium. In comparison, it takes the regular diffusion group less time to reach a lower population density equilibrium that is lower than the initial population density.
In the case where both species survive, under all \textit{1D}, \textit{2D(a)}, and \textit{2D(b)}  cases, small diffusion groups take more time to reach a higher equilibrium. In comparison, the regular diffusion groups take less time reaching a lower equilibrium.

In Figures \ref{cc-2d} and \ref{cc-2ds}, we present solution for 2d problems at final time in whole domain $\Omega$ for regular diffusion $\varepsilon = D$   and  small diffusion $\varepsilon = D/10$, respectively. Similar to Figures \ref{cc-1}, the effect of boundary constraint (set at zero) is more severe in regular diffusion groups. If the diffusion is regular (Figure \ref{cc-2d}), the final population is only dense in the middle, and the area of zero population density is wider at the boundaries. If the diffusion is small (Figure \ref{cc-2ds}), the final population is more spread across the whole domain.
%The possible explanation for this situation might be that if the diffusion rate is larger, species all over the domain can move farther to the center or middle of the domain. At the same time, not all individuals can move to the center of the domain when the diffusion is small.

\subsection{Three-species competition}

For the three-species competition model we consider three cases:
\begin{itemize}
\item Case 1 (one species survive)
\[
D = (0.078, 0.087, 0.012), \quad 
r = (0.050, 0.087, 0.041), \quad
\alpha =
\begin{pmatrix}
0.0 & 0.048 & 0.067\\ 
0.051 & 0.0 & 0.094\\
0.028 & 0.041 & 0.0
\end{pmatrix}
\]
\item Case 2 (both species survive)
\[
D = (0.022, 0.021, 0.063), \quad 
r = (0.086, 0.091, 0.066), \quad
\alpha =
\begin{pmatrix}
0.0 & 0.031 & 0.045\\ 
0.051 & 0.0 & 0.019\\
0.058 & 0.085 & 0.0
\end{pmatrix}
\]
\item Case 3 (all species survive)
\[
D = (0.031, 0.027, 0.026), \quad 
r = (0.098 0.095, 0.078), \quad
\alpha =
\begin{pmatrix}
0.0 & 0.055 & 0.057\\ 
0.095 & 0.0 & 0.031\\
0.070 & 0.022 & 0.0
\end{pmatrix}
\]
\end{itemize}

%\FloatBarrier
\begin{figure}[h!]
\centering
\begin{subfigure}{0.32\textwidth}
\centering
1D\\
\includegraphics[width=1\linewidth]{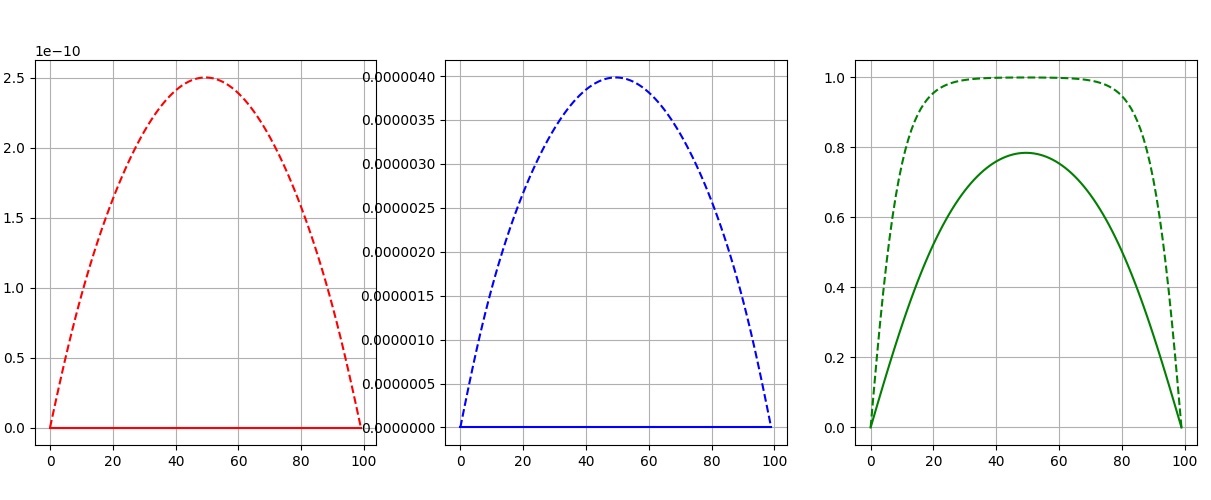}
\end{subfigure}
\begin{subfigure}{0.32\textwidth}
\centering
2D (a)\\
\includegraphics[width=1\linewidth]{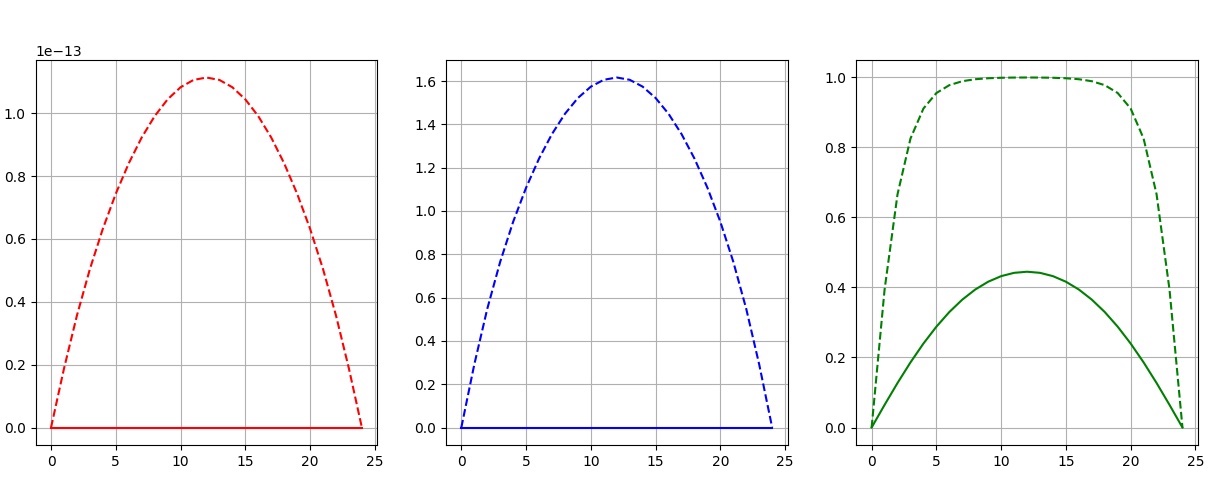}
\end{subfigure}
\begin{subfigure}{0.32\textwidth}
\centering
2D (b)\\
\includegraphics[width=1\linewidth]{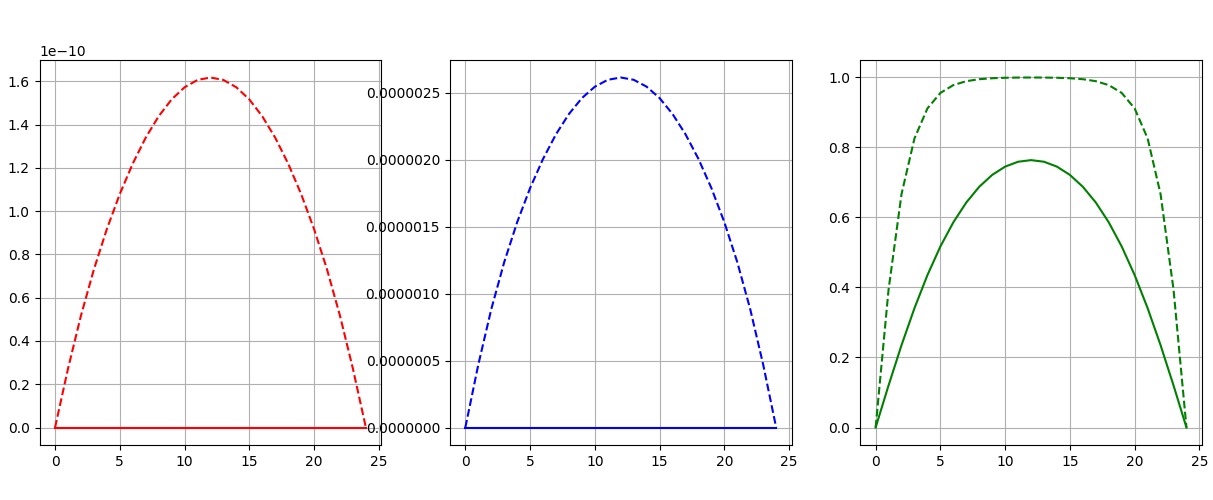}
\end{subfigure}\\
Case 1 (one species survive)
\\
\begin{subfigure}{0.32\textwidth}
\centering
\includegraphics[width=1\linewidth]{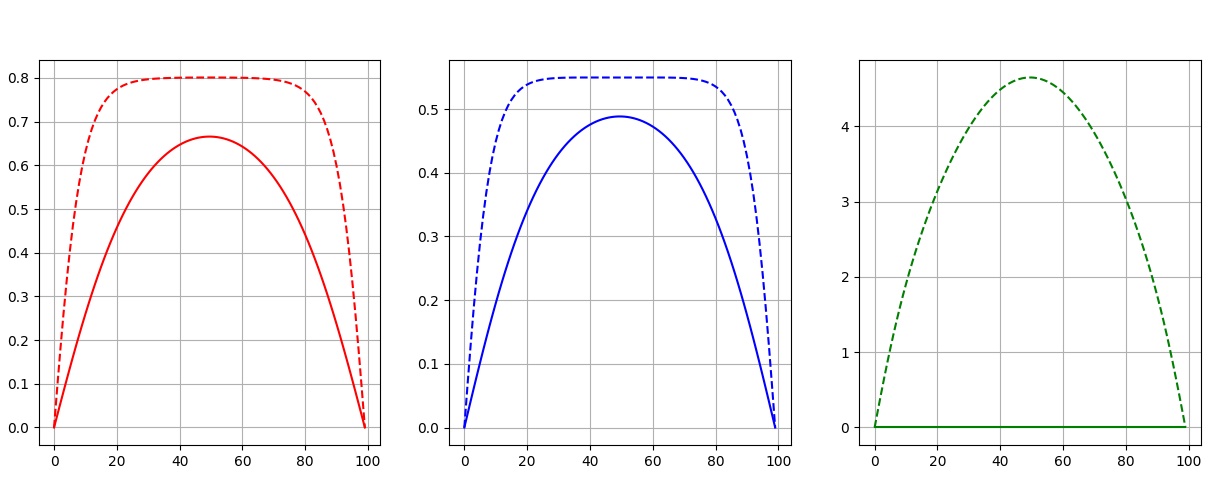}
\end{subfigure}
\begin{subfigure}{0.32\textwidth}
\centering
\includegraphics[width=1\linewidth]{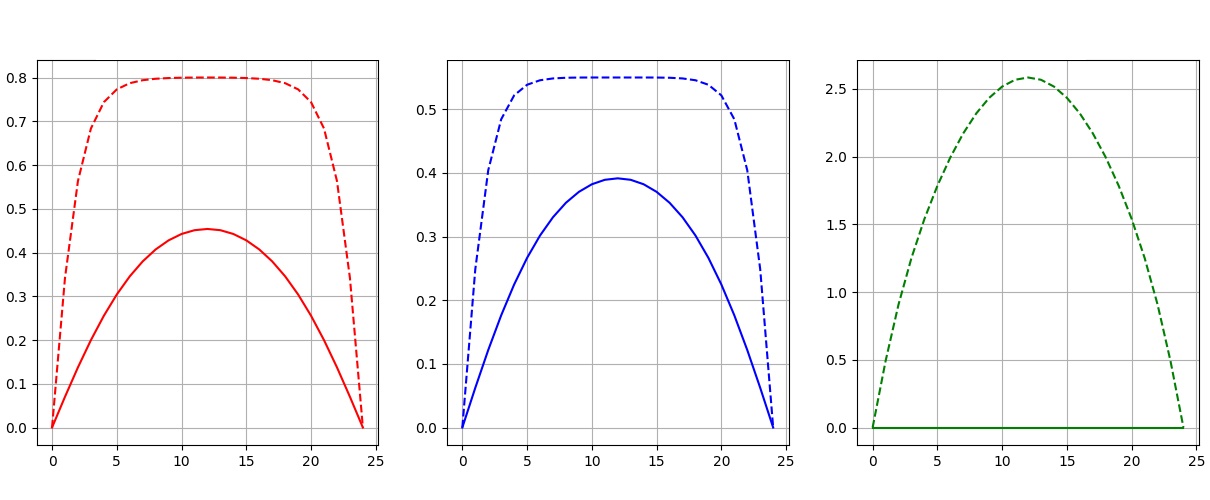}
\end{subfigure}
\begin{subfigure}{0.32\textwidth}
\centering
\includegraphics[width=1\linewidth]{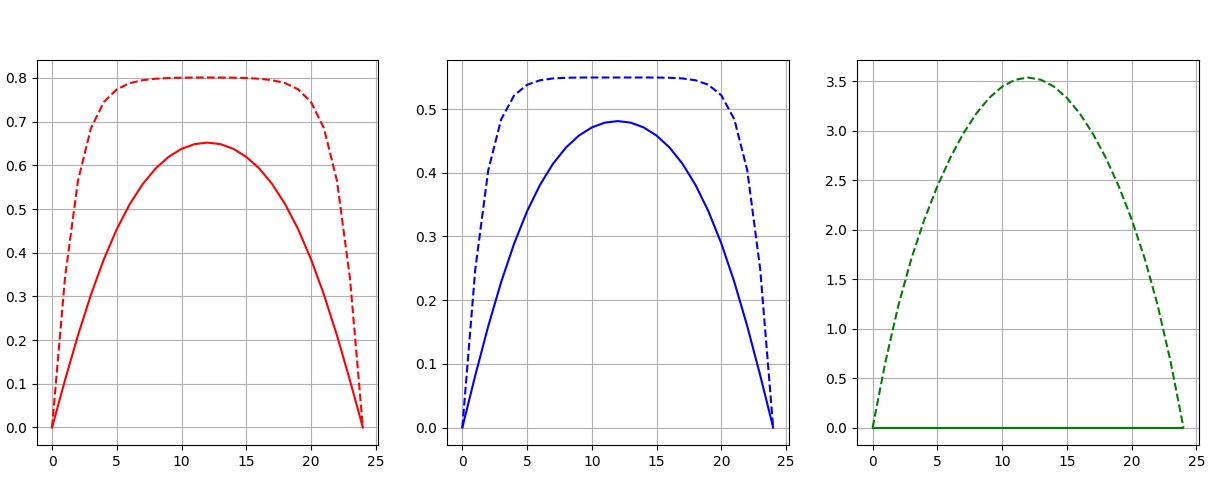}
\end{subfigure}\\
Case 2 (two species survive)
\\
\begin{subfigure}{0.32\textwidth}
\centering
\includegraphics[width=1\linewidth]{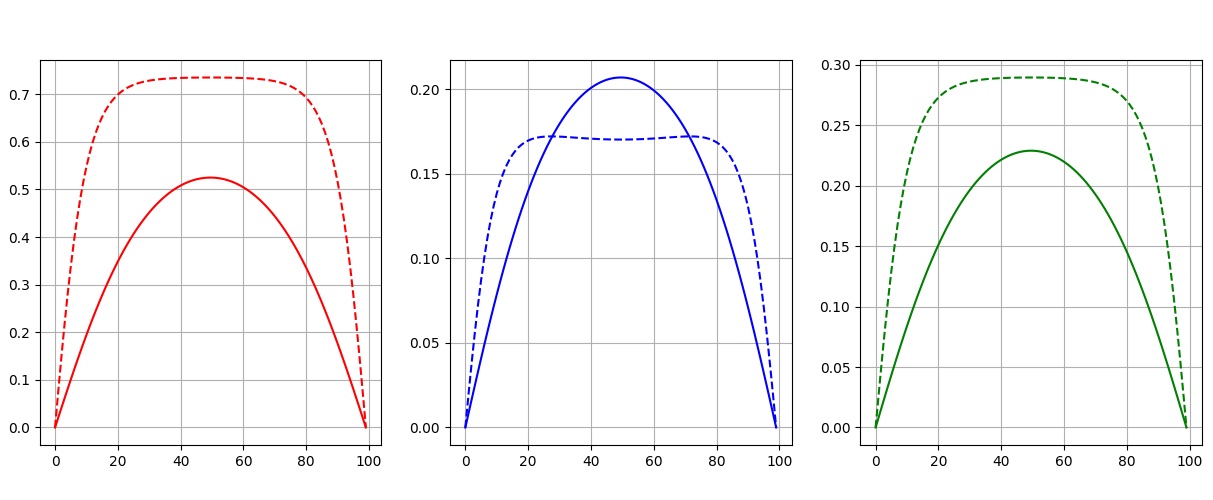}
\end{subfigure}
\begin{subfigure}{0.32\textwidth}
\centering
\includegraphics[width=1\linewidth]{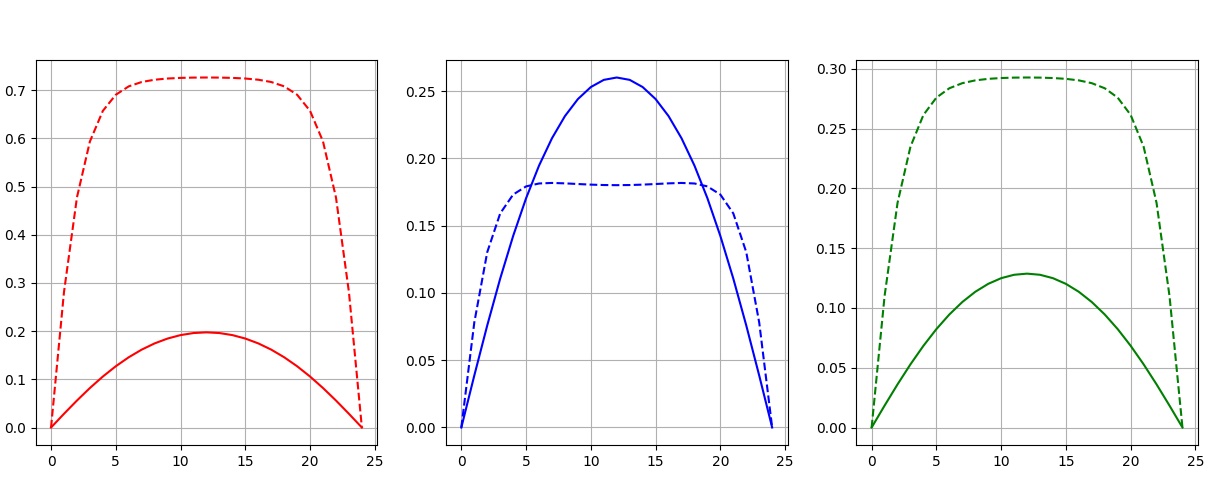}
\end{subfigure}
\begin{subfigure}{0.32\textwidth}
\centering
\includegraphics[width=1\linewidth]{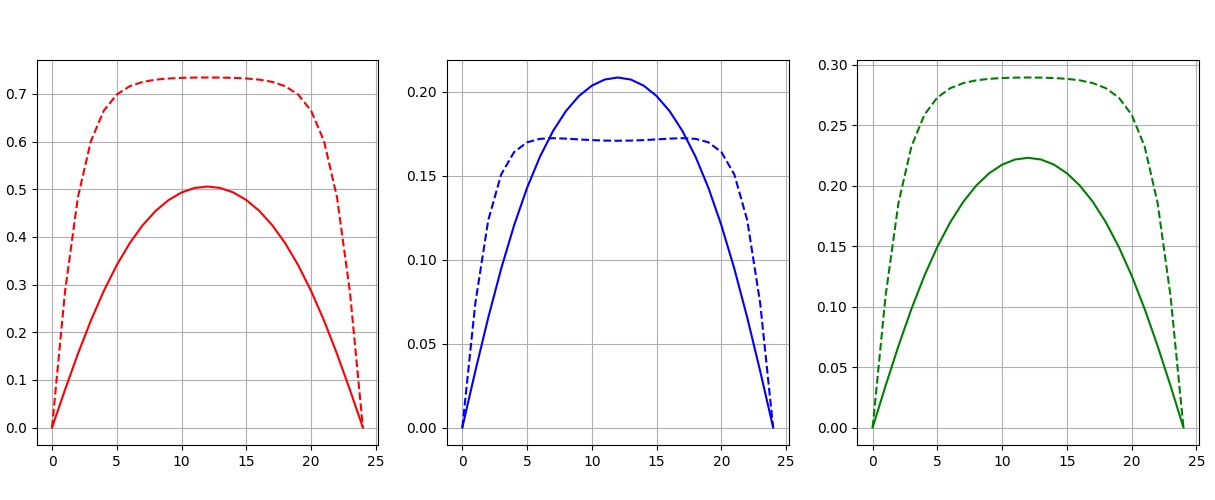}
\end{subfigure}\\
Case 3 (three species survive)
\caption{Solution at final time for regular diffusion $\varepsilon =  D$ (solid line) and small diffusion $\varepsilon = D/10$ (dashed line). 
Red color: first species. Blue color: second species. Green color: third species }
\label{ccc-1}
\end{figure}
%\FloatBarrier

%\FloatBarrier
\begin{figure}[h!]
\centering
\begin{subfigure}{0.32\textwidth}
\centering
1D\\
\includegraphics[width=0.8\linewidth]{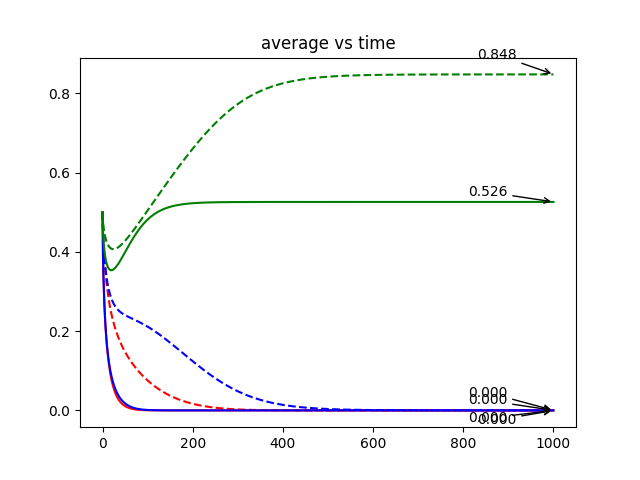}
\end{subfigure}
\begin{subfigure}{0.32\textwidth}
\centering
2D (a)\\
\includegraphics[width=0.8\linewidth]{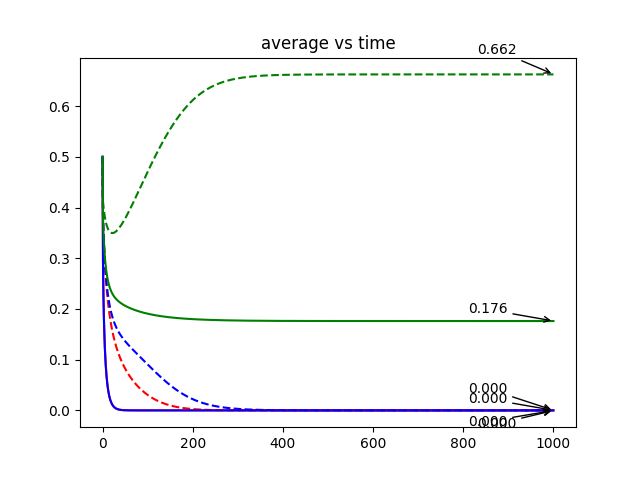}
\end{subfigure}
\begin{subfigure}{0.32\textwidth}
\centering
2D (b)\\
\includegraphics[width=0.8\linewidth]{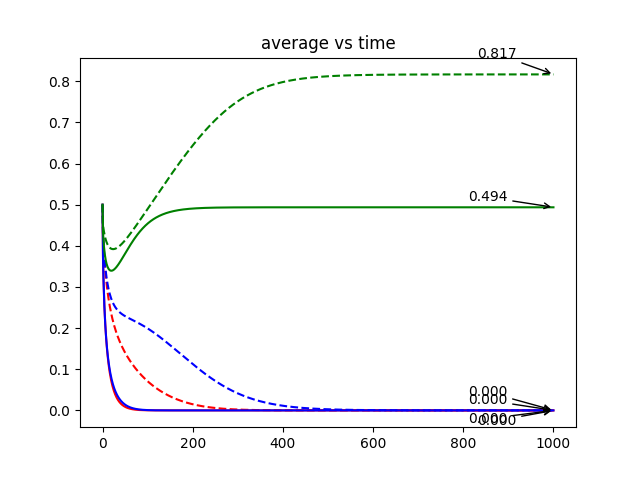}
\end{subfigure}\\
Case 1 (one species survive)
\\
\begin{subfigure}{0.32\textwidth}
\centering
\includegraphics[width=0.8\linewidth]{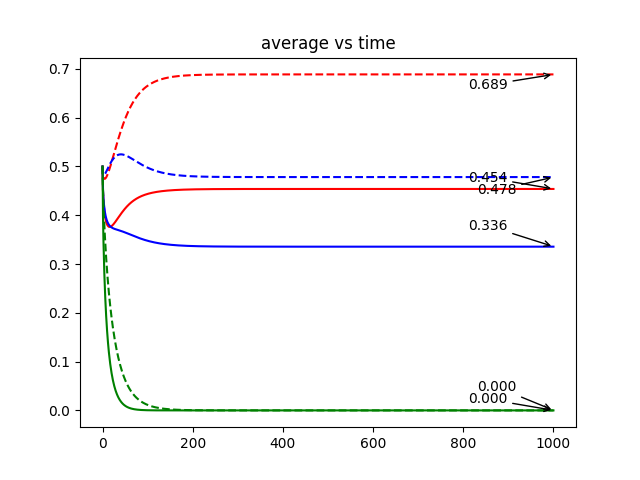}
\end{subfigure}
\begin{subfigure}{0.32\textwidth}
\centering
\includegraphics[width=0.8\linewidth]{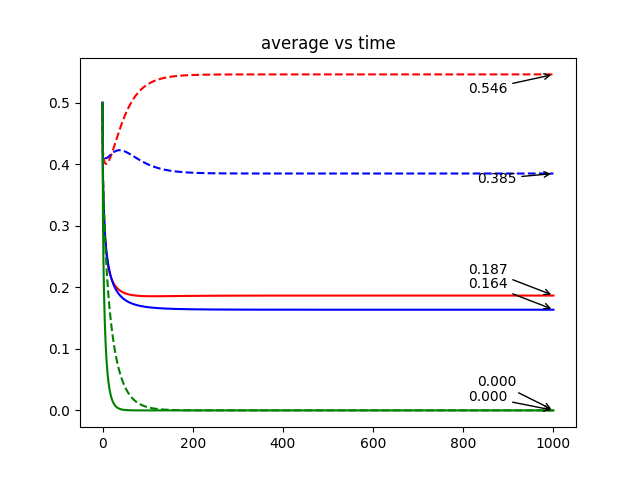}
\end{subfigure}
\begin{subfigure}{0.32\textwidth}
\centering
\includegraphics[width=0.8\linewidth]{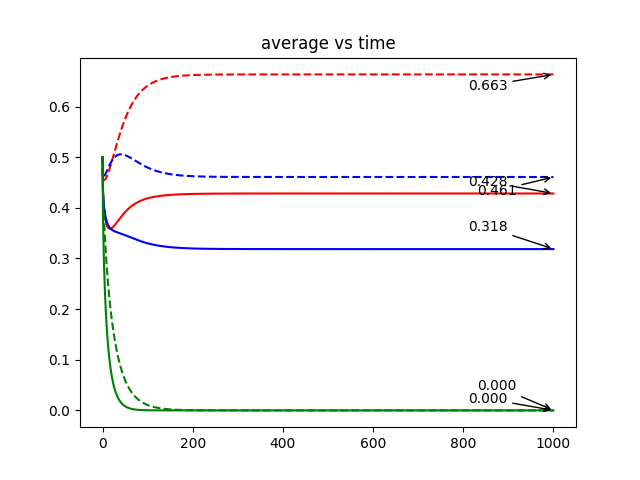}
\end{subfigure}\\
Case 2 (two species survive)
\\
\begin{subfigure}{0.32\textwidth}
\centering
\includegraphics[width=0.8\linewidth]{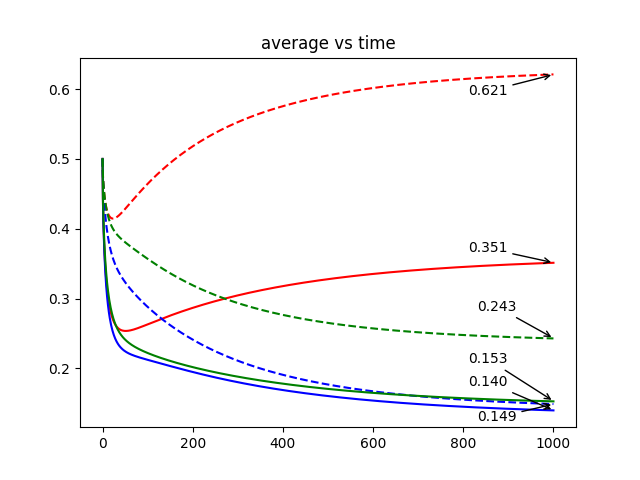}
\end{subfigure}
\begin{subfigure}{0.32\textwidth}
\centering
\includegraphics[width=0.8\linewidth]{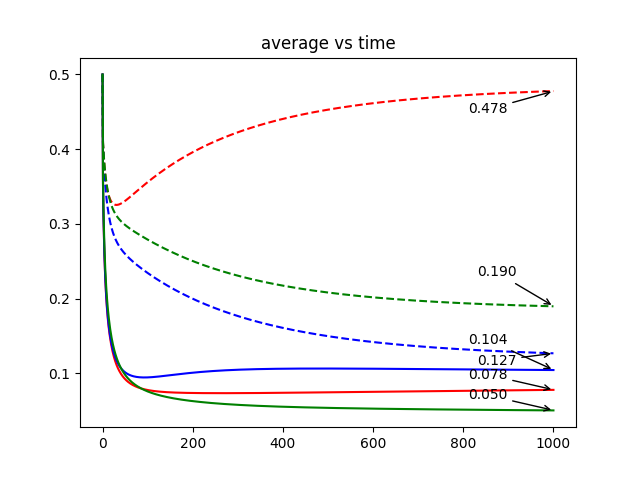}
\end{subfigure}
\begin{subfigure}{0.32\textwidth}
\centering
\includegraphics[width=0.8\linewidth]{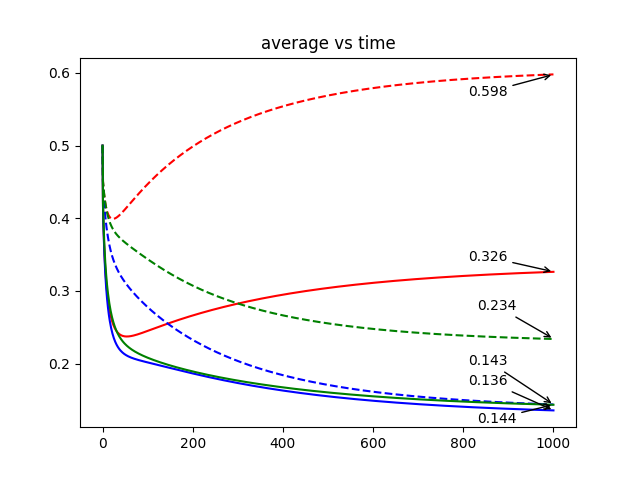}
\end{subfigure}\\
Case 3 (three species survive)
\caption{Solution at final time for regular diffusion $\varepsilon =  D$ (solid line) and small diffusion $\varepsilon = D/10$ (dashed line). 
Red color: first species. Blue color: second species. Green color: third species }
\label{ccc-1t}
\end{figure}
%\FloatBarrier

In Figures \ref{ccc-1}, we plot the solution for one and two-dimensional formulations at the final time. Similarly to the previous results, we plot the solution over the middle line ($y = 0.5$) for the 2d formulation. We observe that the effect of boundary conditions is generally more severe in regular diffusion groups than in small diffusion groups. %The zero final population density area is wider for the regular diffusion group. 
With the same survival status (one species survives, or both species survive) and spatial boundary conditions (\textit{1D}, \textit{2D(a)}, or \textit{2D(b)}), the regular diffusion group mostly arrive at a lower population density equilibrium, compare to the small diffusion group.  The species can reach a higher population density equilibrium when the diffusion is smaller. 
However, in contrast to the two species system, where the boundary effect holds for all species in every survival status scenario, there is a special case in the three species system. When all three species survive (case 3 in Figures \ref{ccc-1}), one of the species reaches a higher final population density when the diffusion is regular instead of when the diffusion is small. More investigation of the diffusion to the final state will be presented in the next part of the results. 

In Figures \ref{ccc-1t}, we present the dynamic of the solution average over time for three species system. We observed that the time it takes to reach equilibrium is different among different diffusion conditions (regular diffusion $D$, or small diffusion $D/10$), survival status (one species survives, or both species survive) and spatial boundary conditions (\textit{1D}, \textit{2D(a)}, or \textit{2D(b)}). In general, similar to two species system, \textit{1D} and \textit{2D(b)} give similar solutions, while \textit{2D(a)} gives different solution.  Furthermore in contrast to the two species system, where the population density of some species increases from the very beginning, under the three species system, all species decrease in population density from the initial density of 0.5 at first, and then either increase or decrease to an equilibrium.
Under all survival statuses in the three species system, small diffusion groups spend more time to reach a higher (for the survived species) equilibrium, while the regular diffusion groups spend less time to reach a lower equilibrium.

\subsection{Effect of the diffusion}

Next, we consider the influence of diffusion on the equilibrium state. We set all other parameters (growth rate "$r$" and competition efficiency "$\alpha$") the same and examine the diffusion rate.  
We perform 1000 simulations for each case with random diffusion coefficients   
\[
0.01 < \varepsilon_i< 0.1.
\]
We simulate with  fixed initial conditions $u_0^{(k)} = 0.5$  till both population reach equilibrium, $|\bar{u}^{(k)} - \check{\bar{u}}^{(k)}| < \epsilon$ with $\epsilon = 10^{-5}$ for each $k$.  

%\FloatBarrier
\begin{figure}[h!]
\centering
\begin{subfigure}{0.32\textwidth}
\centering
1D\\
\includegraphics[width=0.43\linewidth]{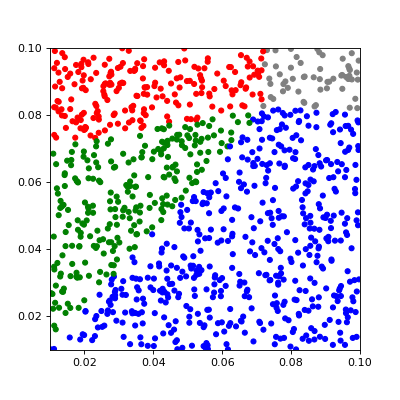}
\includegraphics[width=0.54\linewidth]{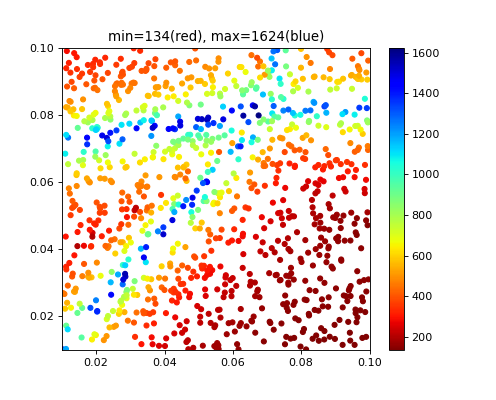}
\end{subfigure}
\begin{subfigure}{0.32\textwidth}
\centering
2D (a)\\
\includegraphics[width=0.43\linewidth]{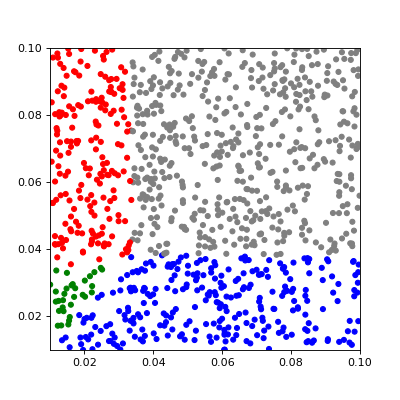}
\includegraphics[width=0.54\linewidth]{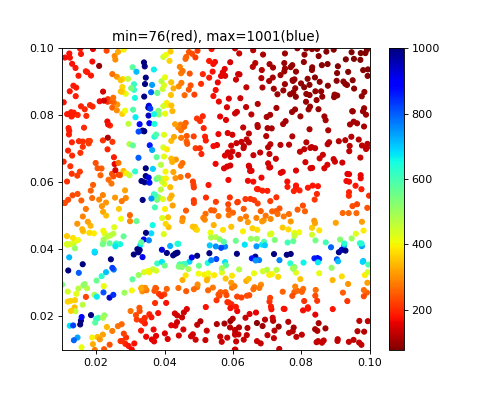}
\end{subfigure}
\begin{subfigure}{0.32\textwidth}
\centering
2D (b)\\
\includegraphics[width=0.43\linewidth]{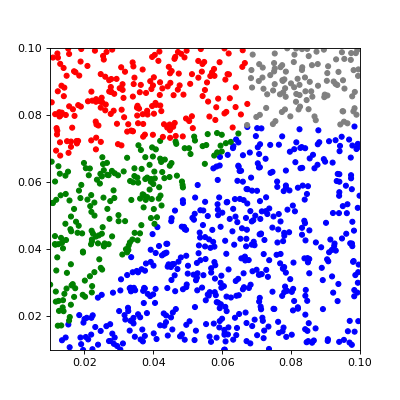}
\includegraphics[width=0.54\linewidth]{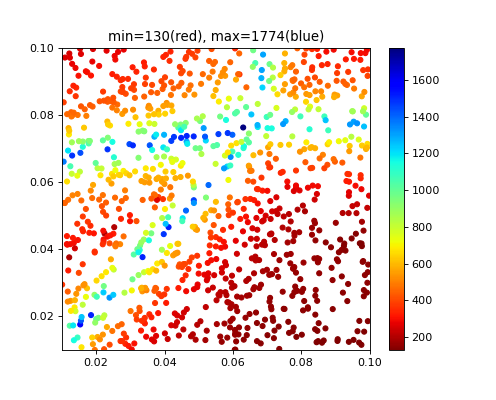}
\end{subfigure}\\
Case 1 (one species survive)
\\
\begin{subfigure}{0.32\textwidth}
\centering
\includegraphics[width=0.43\linewidth]{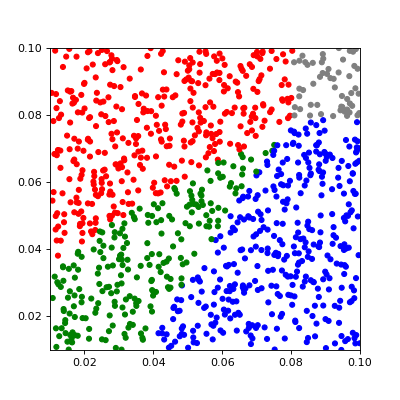}
\includegraphics[width=0.54\linewidth]{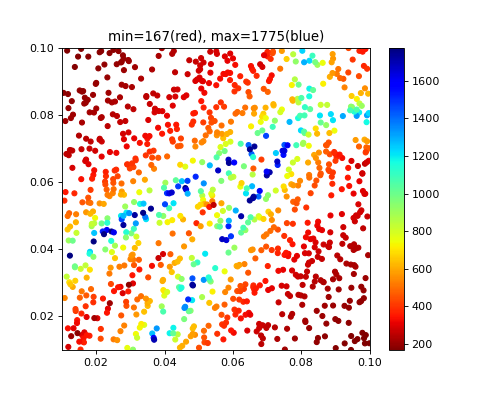}
\end{subfigure}
\begin{subfigure}{0.32\textwidth}
\centering
\includegraphics[width=0.43\linewidth]{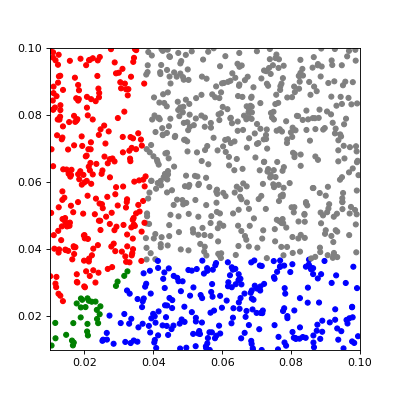}
\includegraphics[width=0.54\linewidth]{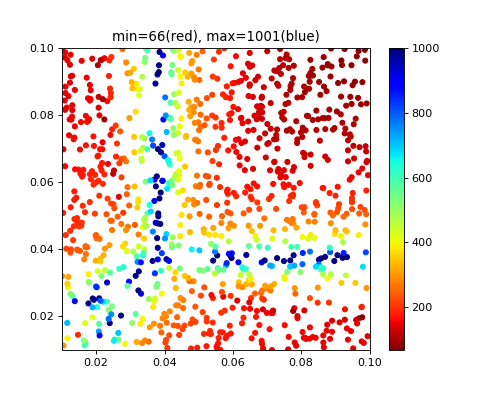}
\end{subfigure}
\begin{subfigure}{0.32\textwidth}
\centering
\includegraphics[width=0.43\linewidth]{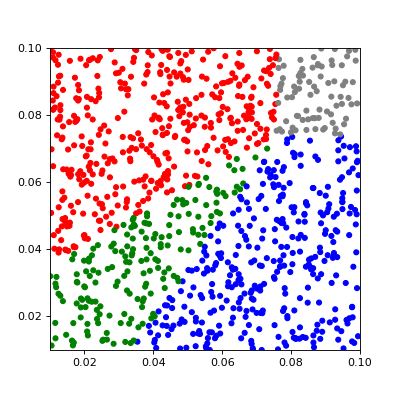}
\includegraphics[width=0.54\linewidth]{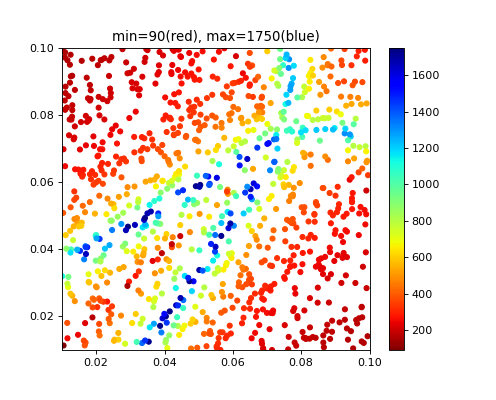}
\end{subfigure}\\
Case 2 (two species survive)
\caption{Numerical results for random diffusion.  Two-species model. 
First picture: 00 (grey),  01 (blue),  10 (red) and 11 (green) groups.
Second picture: number of time steps to reach the equilibrium}
\label{cc-diffusion}
\end{figure}
%\FloatBarrier

%\FloatBarrier
\begin{figure}[h!]
\centering
\begin{subfigure}{0.32\textwidth}
\centering
1D\\
\includegraphics[width=0.48\linewidth]{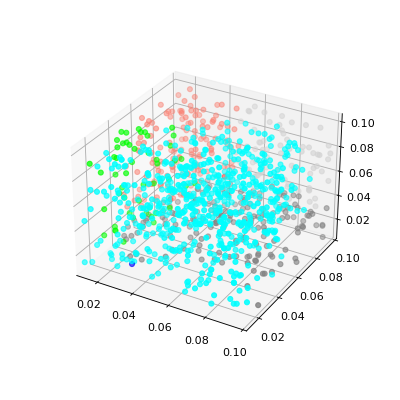}
\includegraphics[width=0.48\linewidth]{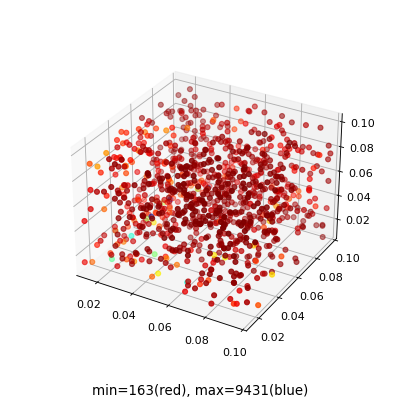}
\end{subfigure}
\begin{subfigure}{0.32\textwidth}
\centering
2D (a)\\
\includegraphics[width=0.48\linewidth]{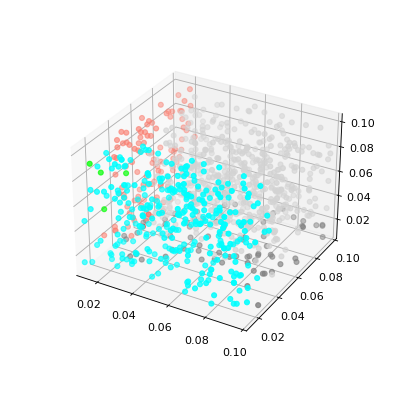}
\includegraphics[width=0.48\linewidth]{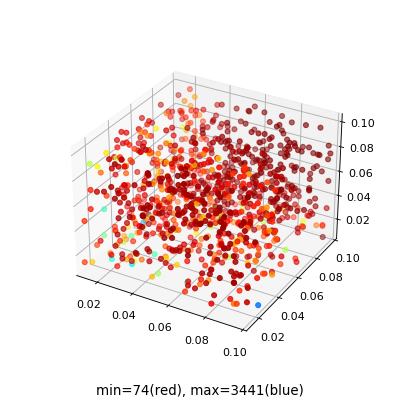}
\end{subfigure}
\begin{subfigure}{0.32\textwidth}
\centering
2D (b)\\
\includegraphics[width=0.48\linewidth]{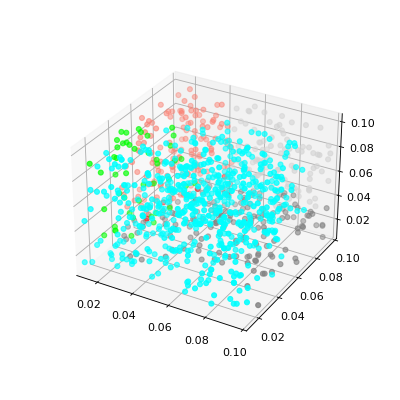}
\includegraphics[width=0.48\linewidth]{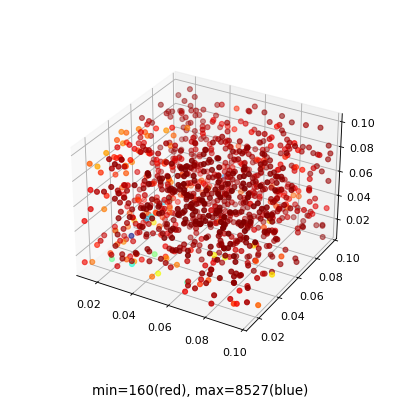}
\end{subfigure}\\
Case 1 (one species survive)
\\
\begin{subfigure}{0.32\textwidth}
\centering
\includegraphics[width=0.48\linewidth]{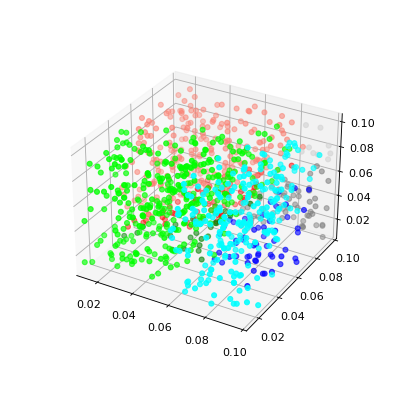}
\includegraphics[width=0.48\linewidth]{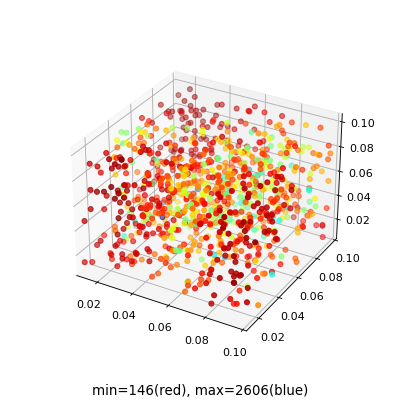}
\end{subfigure}
\begin{subfigure}{0.32\textwidth}
\centering
\includegraphics[width=0.48\linewidth]{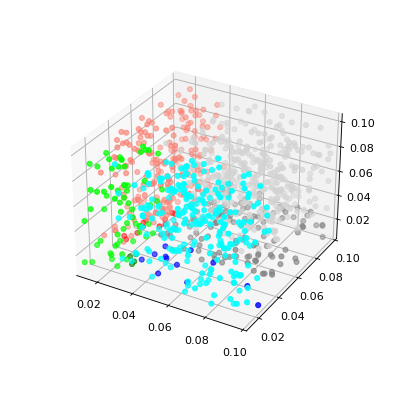}
\includegraphics[width=0.48\linewidth]{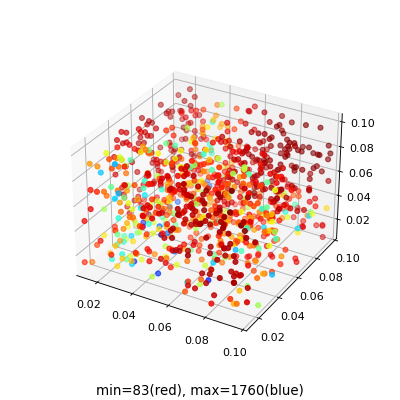}
\end{subfigure}
\begin{subfigure}{0.32\textwidth}
\centering
\includegraphics[width=0.48\linewidth]{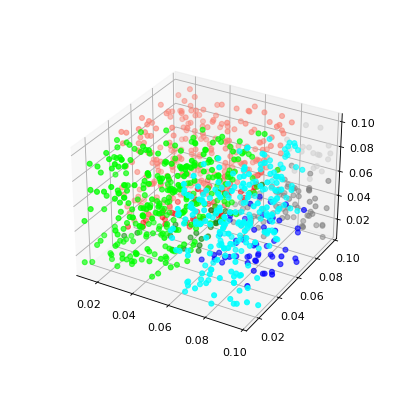}
\includegraphics[width=0.48\linewidth]{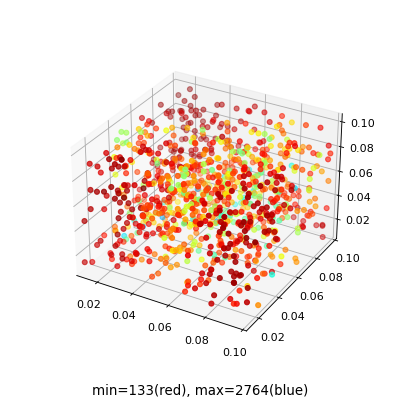}
\end{subfigure}\\
Case 2 (two species survive)
\\
\begin{subfigure}{0.32\textwidth}
\centering
\includegraphics[width=0.48\linewidth]{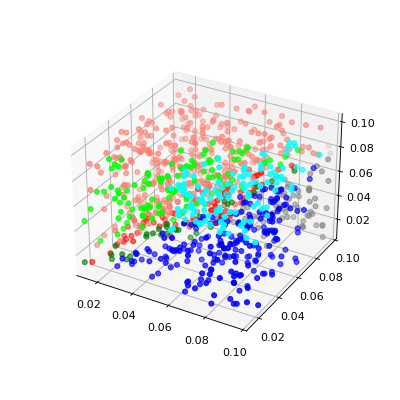}
\includegraphics[width=0.48\linewidth]{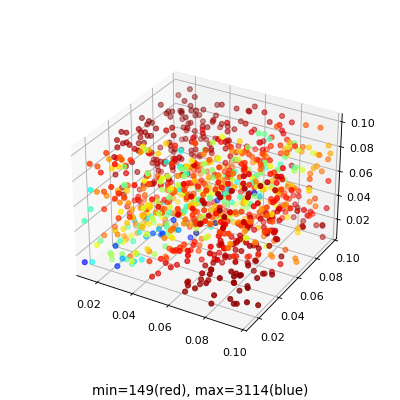}
\end{subfigure}
\begin{subfigure}{0.32\textwidth}
\centering
\includegraphics[width=0.48\linewidth]{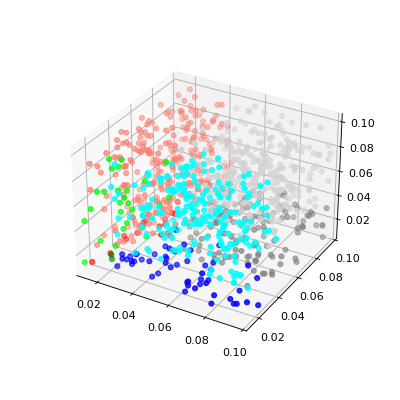}
\includegraphics[width=0.48\linewidth]{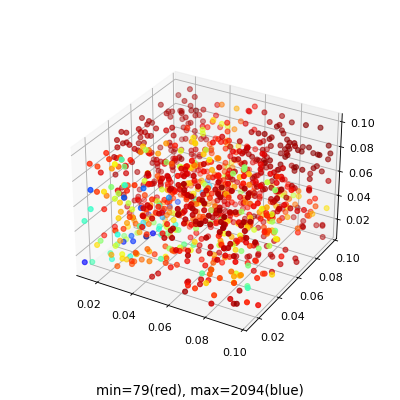}
\end{subfigure}
\begin{subfigure}{0.32\textwidth}
\centering
\includegraphics[width=0.48\linewidth]{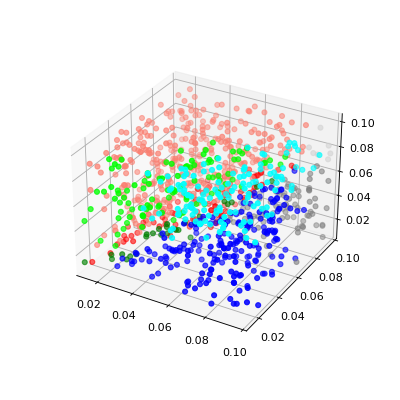}
\includegraphics[width=0.48\linewidth]{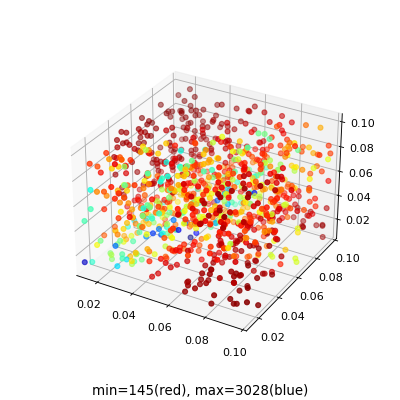}
\end{subfigure}\\
Case 3 (three species survive)
\caption{Numerical results for random diffusion.  Three-species model. 
First picture: 000 (light grey),  001 (grey), 010 (cyan),  011 (blue),  100 (salmon),  101 (red),  110 (lime) and 111 (green) groups.
Second picture: number of time steps to reach the equilibrium}
\label{ccc-diffusion}
\end{figure}
%\FloatBarrier

First, we consider two-species competition models in one and two-dimensional formulations.   
Groups are represented by which species survive.  
In the two-species competition model we have: 
10 - first species survive,  
01 - second species survived, 
11 - both species survived and  
00 - no one survived.

In Figure \ref{cc-diffusion}, we plot scatter plots of the diffusion rate of two species and colored survival status at the final time for one and two-dimensional problems, as well as corresponding time steps to reach equilibrium. We observed that diffusion rate have a huge impacts to the final survival status of species, as well as the time step needed to reach equilibrium. 

Different combination of diffusion rates of two species leads to different final survival status. Under \textit{1D} or \textit{2D(b)} spatial boundary condition, when the diffusion of both species are above 0.07, no species will survive. When one of the species' diffusion rates is smaller than 0.07, the species with a larger diffusion rate will survive. However, if two species have approximately the same diffusion rates and are both less than 0.07, both species will survive. In comparison, under \textit{2D(a)} spatial boundary condition, the threshold is reduced from 0.07 to around 0.04. 
Different combination of diffusion rates of two species also leads to different time steps needed to reach equilibrium. When the combination of both species' diffusion rates is on the border of different survival groups, it takes more time steps to reach equilibrium, while it takes fewer time steps to reach equilibrium when the combination of diffusion rates of both species lies inside the survival group.

Next, we consider three-species competition models, where we have eight groups: 
100 - first species survive,  
010 - second species survived, 
001 - third species survived, 
011 - second and third survived, 
101 - first and third survived, 
110 - first and  second survived, 
111 - all species survived and  
000 - no one survived.

In Figure \ref{ccc-diffusion} we plot scatter plot of diffusion rate of three species and colored survival status at final time for one and two-dimensional spatial boundary conditions, as well as corresponding time steps to reach equilibrium. When we expand from two species to three species in the system, we observed that \textit{1D} or \textit{2D(b)} spatial boundary condition still gives similar result, while \textit{2D(a)} gives different result. In general, when the diffusion rates of all three species are small, three species can all survive, when all species have high diffusion rates, no species can survive. When not all three species have low diffusion rates, the species with higher diffusion rate can survive.
Similarly,   the different combination of diffusion rates of the three species also leads to different time steps needed to reach equilibrium and maximum located on the border of different survival groups.

\subsection{Effect of the initial conditions}

%Next, we consider influence of the initial conditions to the final solution (equilibrium state).  We consider two and three-species competition models in one and two-dimensional formulations  with fixed diffusion rate $\varepsilon = D$.  

Previously we set the initial condition for both species to be 0.5 and then get a general sense of the final equilibrium difference.  Next we remove the control for the initial condition ($0.01 < u_0^{(k)} < 0.99$), and examine how changes of the initial population can impact the final equilibrium population density and time to reach it. With the birth rate, competition rate, and diffusion rate fixed, we change the initial population density for both species.

%\FloatBarrier
\begin{figure}[h!]
\centering
\begin{subfigure}{0.32\textwidth}
\centering
1D\\
\includegraphics[width=0.45\linewidth]{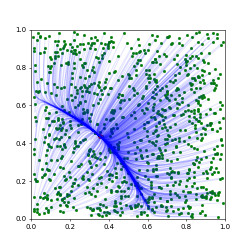}
\includegraphics[width=0.52\linewidth]{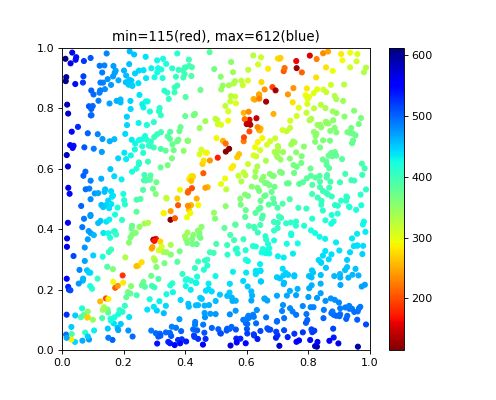}
\end{subfigure}
\begin{subfigure}{0.32\textwidth}
\centering
2D (a)\\
\includegraphics[width=0.45\linewidth]{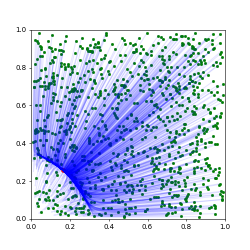}
\includegraphics[width=0.52\linewidth]{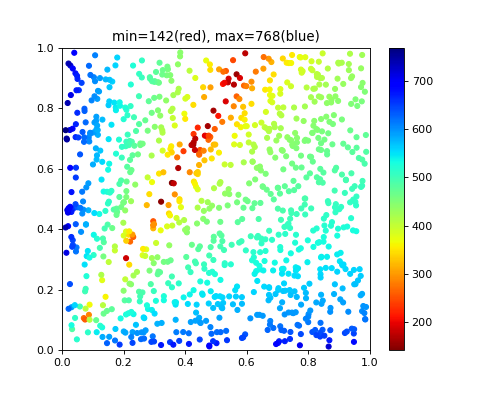}
\end{subfigure}
\begin{subfigure}{0.32\textwidth}
\centering
2D (b)\\
\includegraphics[width=0.45\linewidth]{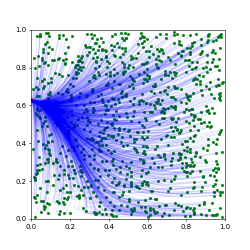}
\includegraphics[width=0.52\linewidth]{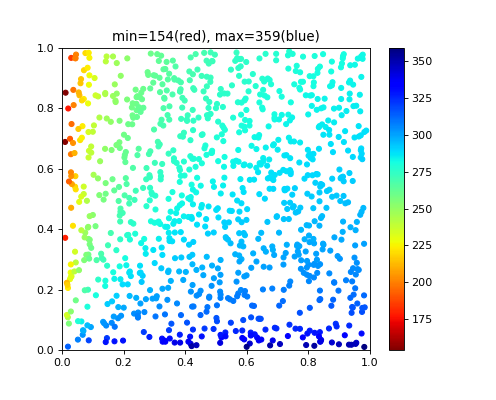}
\end{subfigure}
\\
Case 1 (one species survive)
\\
\begin{subfigure}{0.32\textwidth}
\centering
\includegraphics[width=0.45\linewidth]{new_random_IC/sim_ic_3-a}
\includegraphics[width=0.52\linewidth]{new_random_IC/sim_ic_3-b}
\end{subfigure}
\begin{subfigure}{0.32\textwidth}
\centering
\includegraphics[width=0.45\linewidth]{new_random_IC/sim_ic_3_2d_a-a}
\includegraphics[width=0.52\linewidth]{new_random_IC/sim_ic_3_2d_a-b}
\end{subfigure}
\begin{subfigure}{0.32\textwidth}
\centering
\includegraphics[width=0.45\linewidth]{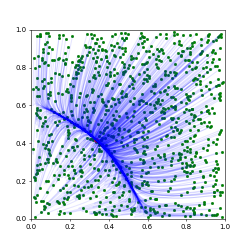}
\includegraphics[width=0.52\linewidth]{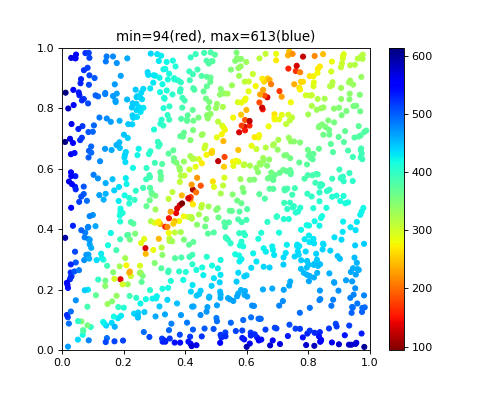}
\end{subfigure}\\
Case 2 (two species survive)
\caption{Numerical results for random initial conditions.  Two-species model. 
First picture: dynamic of the solution average.
Second picture: number of time steps to reach the equilibrium}
\label{cc-ic}
\end{figure}
%\FloatBarrier

%\FloatBarrier
\begin{figure}[h!]
\centering
\begin{subfigure}{0.32\textwidth}
\centering
1D\\
\includegraphics[width=0.48\linewidth]{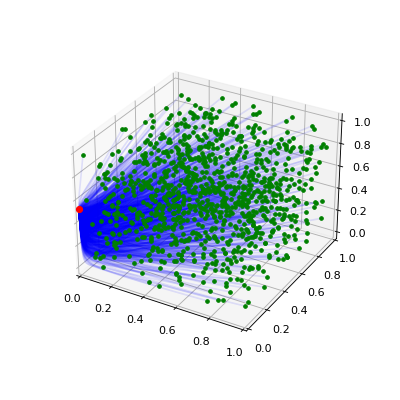}
\includegraphics[width=0.48\linewidth]{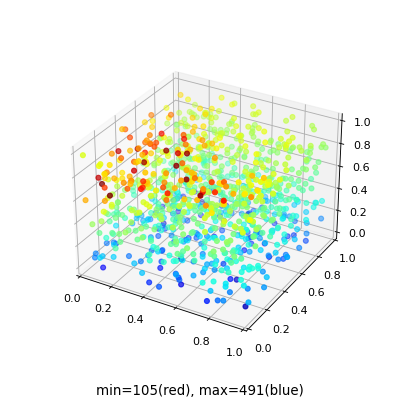}
\end{subfigure}
\begin{subfigure}{0.32\textwidth}
\centering
2D (a)\\
\includegraphics[width=0.48\linewidth]{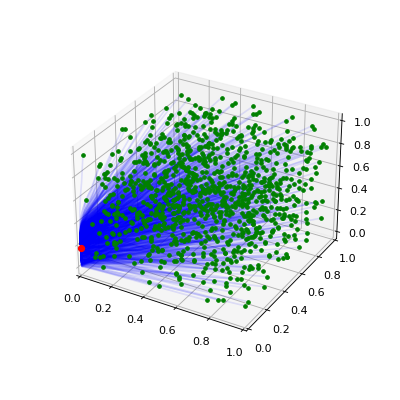}
\includegraphics[width=0.48\linewidth]{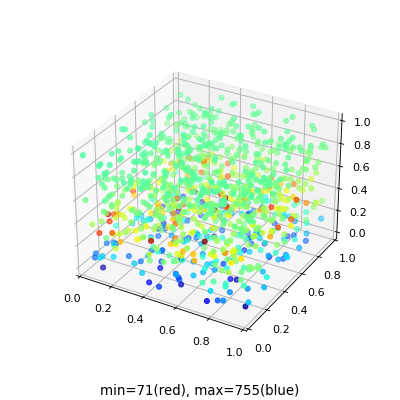}
\end{subfigure}
\begin{subfigure}{0.32\textwidth}
\centering
2D (b)\\
\includegraphics[width=0.48\linewidth]{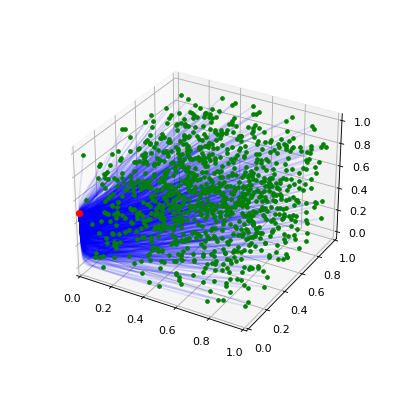}
\includegraphics[width=0.48\linewidth]{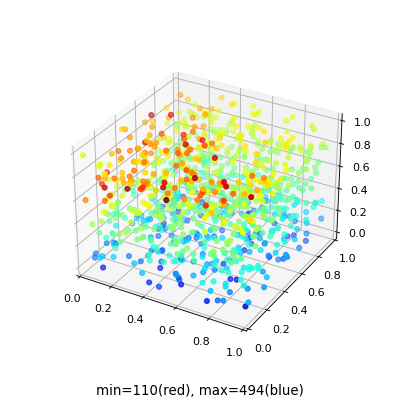}
\end{subfigure}\\
Case 1 (one species survive)
\\
\begin{subfigure}{0.32\textwidth}
\centering
\includegraphics[width=0.48\linewidth]{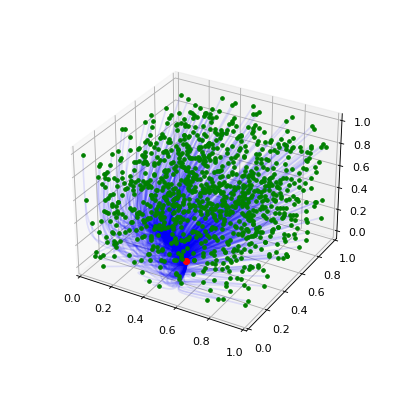}
\includegraphics[width=0.48\linewidth]{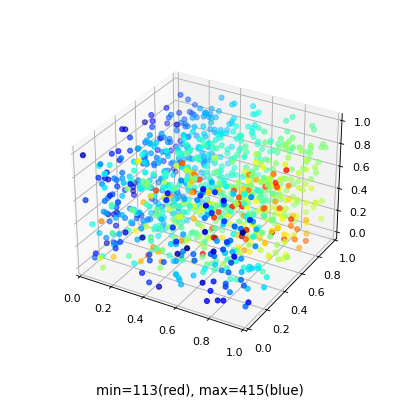}
\end{subfigure}
\begin{subfigure}{0.32\textwidth}
\centering
\includegraphics[width=0.48\linewidth]{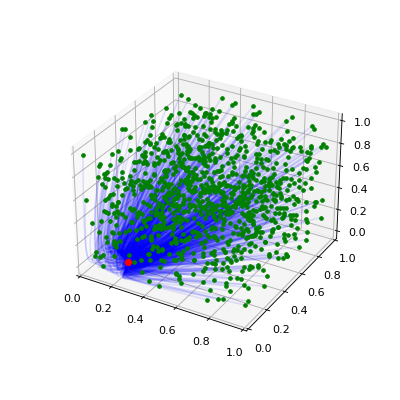}
\includegraphics[width=0.48\linewidth]{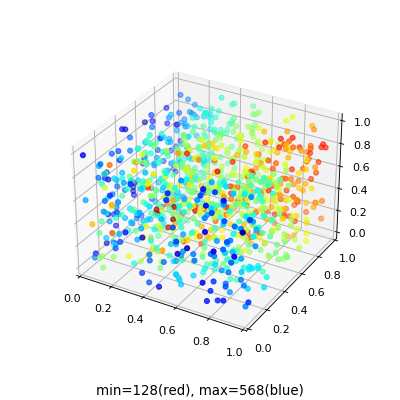}
\end{subfigure}
\begin{subfigure}{0.32\textwidth}
\centering
\includegraphics[width=0.48\linewidth]{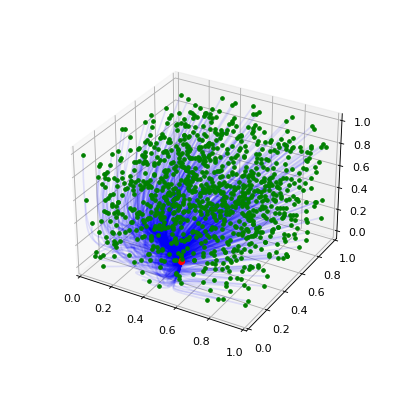}
\includegraphics[width=0.48\linewidth]{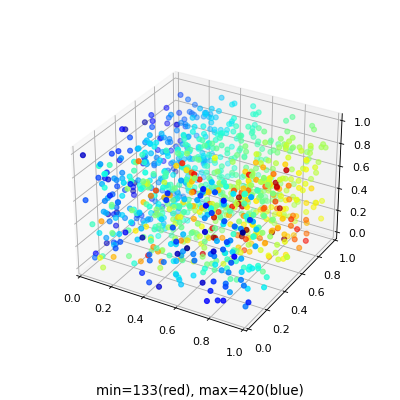}
\end{subfigure}\\
Case 2 (two species survive)
\\
\begin{subfigure}{0.32\textwidth}
\centering
\includegraphics[width=0.48\linewidth]{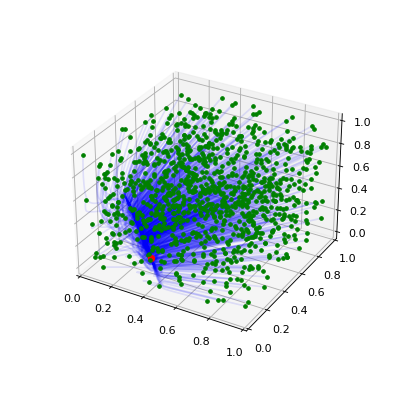}
\includegraphics[width=0.48\linewidth]{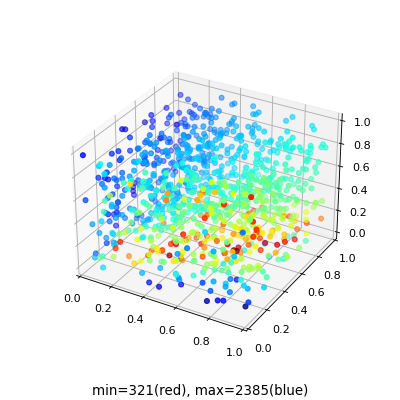}
\end{subfigure}
\begin{subfigure}{0.32\textwidth}
\centering
\includegraphics[width=0.48\linewidth]{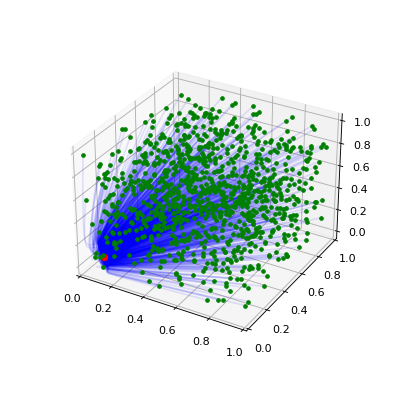}
\includegraphics[width=0.48\linewidth]{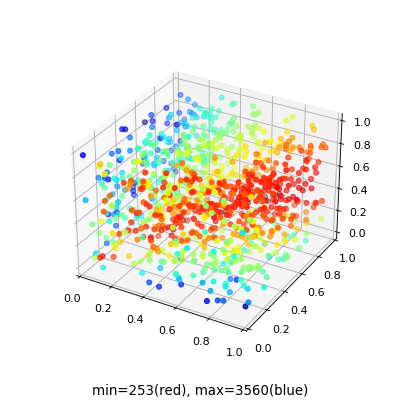}
\end{subfigure}
\begin{subfigure}{0.32\textwidth}
\centering
\includegraphics[width=0.48\linewidth]{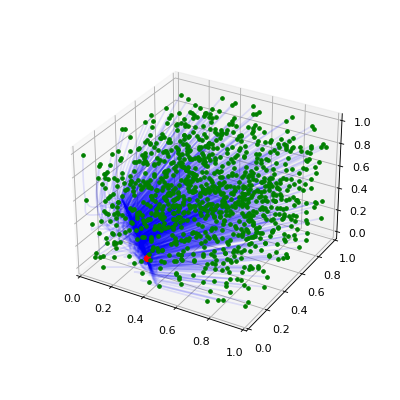}
\includegraphics[width=0.48\linewidth]{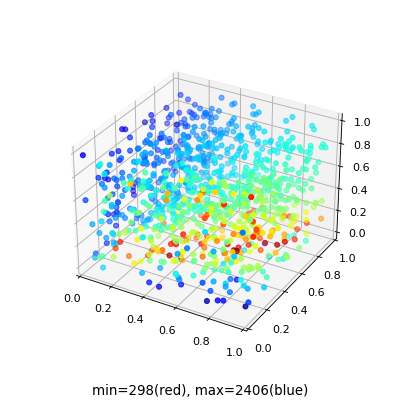}
\end{subfigure}\\
Case 3 (three species survive)
\caption{Numerical results for random initial conditions.  Three-species model. 
First picture: dynamic of the solution average.
Second picture: number of time steps to reach the equilibrium}
\label{ccc-ic}
\end{figure}
%\FloatBarrier

In Figure \ref{cc-ic} we show a scatter plot of the initial population density of two species and use a blue line to represent a dynamic from initial condition points (green color) to final equilibrium points of each simulation for one and two-dimensional spatial boundary conditions, as well as corresponding time steps to reach equilibrium.
We observed that disregarding the value of the initial population of both species, the final equilibrium population will rest at a point of focus. The farther the combination point of initial population density is from the final equilibrium combination point, the more time steps it would take to reach the equilibrium.
We also observed that the farther the initial population density point is from the final equilibrium population density point, the more time steps it needs to reach the equilibrium.

In Figure \ref{ccc-ic}, we plot scatter plots of the initial population density of three species and use a blue line to connect initial condition points to the final equilibrium points of each simulation and corresponding time steps to reach equilibrium.
When we expand from two species to three species in the system, we observed that \textit{1D} or \textit{2D(b)} spatial boundary condition still gives a similar result. In contrast, \textit{2D(a)} gives a different result. In general, regardless of the values of the initial condition of population density of three species, the final equilibrium of population density for all three species in the system will rest at a final focus point. The farther a combination of the initial condition of population density is from the final equilibrium point, the more time steps it needs to take to reach that equilibrium point.

%\clearpage
\section{Factor Analysis}
Factor analysis method has long been applied to analysis of Population Dynamic \cite{varleyRecentAdvancesInsect1970,albonTemporalChangesKey2000} as well as Ecosystem topic \cite{bernierStructuralRelationshipsVegetation2012,petitgasEcosystemSpatialStructure2018}.
Finally, we present factor analysis for the given model. 

Previous literature about factor analysis suggests that as sample size increases, the standard error in factor loadings across repeated samples will decrease \cite{maccallumSampleSizeFactor1999}. Therefore, we perform 10,000 simulations for each case with random parameters to satisfy the sample size requirement.

We consider one and two dimensional test problems and simulate the system with random  coefficients $\varepsilon_k$,  $r_k$, $\alpha_{kl}$ with $l \neq k$ and initial conditions $u^{(k)}_0 = \const$:
\begin{itemize}
\item \textit{Two-species}: 
$\varepsilon_1$, $\varepsilon_2$,  
$\alpha_{12}$, $\alpha_{21}$ 
and $u_0 = [u_0^{(1)}, u_0^{(2)}]$.
\item \textit{Three-species}: 
$\varepsilon_1$, $\varepsilon_2$, $\varepsilon_3$, 
$\alpha_{12}$, $\alpha_{13}$, $\alpha_{21}$, $\alpha_{23}$, 
$\alpha_{31}$, $\alpha_{32}$ 
and $u_0 = [u_0^{(1)}, u_0^{(2)} u_0^{(3)}]$.
\end{itemize}

We perform 10,000 simulations for each case with random parameters 
\[
0.01 < r_k, D_k, \alpha_{kl}, < 0.1,
\]
and $0.01 < u_0^{(k)} < 0.99$. 
Two cases of the diffusion coefficient is considered: 
regular diffusion ($\varepsilon_k = D_k$) and small diffusion($\varepsilon_k = D_k / 10$). 

%\FloatBarrier
\begin{figure}[h!]
\centering
\begin{subfigure}{0.32\textwidth}
\centering
1D\\
\includegraphics[width=0.99\linewidth]{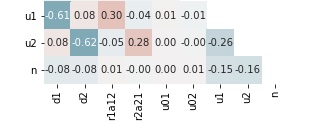}\\
\includegraphics[width=0.99\linewidth]{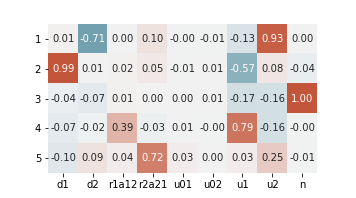}
\end{subfigure}
\begin{subfigure}{0.32\textwidth}
\centering
2D (a)\\
\includegraphics[width=0.99\linewidth]{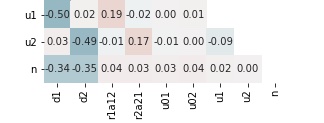}
\includegraphics[width=0.99\linewidth]{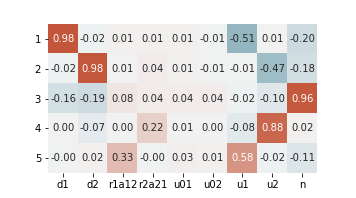}
\end{subfigure}
\begin{subfigure}{0.32\textwidth}
\centering
2D (b)\\
\includegraphics[width=0.99\linewidth]{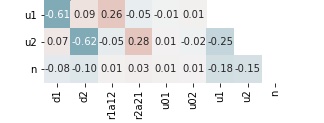}
\includegraphics[width=0.99\linewidth]{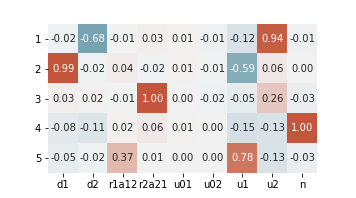}
\end{subfigure}\\
Regular diffusion
\\
\vspace{10pt}
\begin{subfigure}{0.32\textwidth}
\centering
\includegraphics[width=0.99\linewidth]{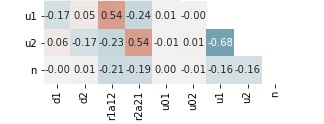}
\includegraphics[width=0.99\linewidth]{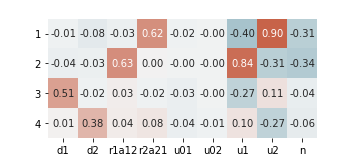}
\end{subfigure}
\begin{subfigure}{0.32\textwidth}
\centering
\includegraphics[width=0.99\linewidth]{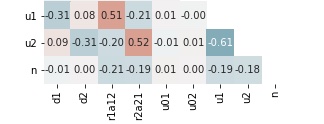}
\includegraphics[width=0.99\linewidth]{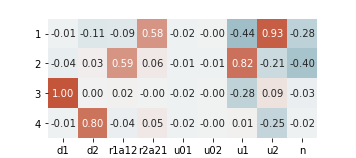}
\end{subfigure}
\begin{subfigure}{0.32\textwidth}
\centering
\includegraphics[width=0.99\linewidth]{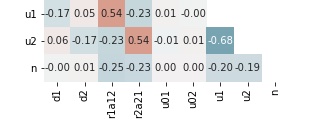}
\includegraphics[width=0.99\linewidth]{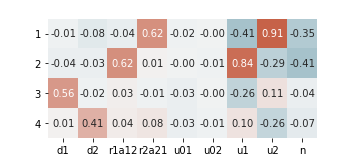}
\end{subfigure}\\
Small diffusion
\caption{Two-species competition model. Correlation matrix (first row) and loading  (second row) }
\label{cc-loading}
\end{figure}
%\FloatBarrier

In investigation we use influence of the following  parameters to the system solution: $\varepsilon_k$ and $\alpha_{kl}/r_k$ for $l \neq k$ and $k=1,2$ for two-species competition and $k=1,2,3$ for three species competition.
In Figure \ref{cc-loading}, we present correlation matrices as well as loading matrix for each factor analysis for all three cases \textit{1D}, \textit{2D(a)} and \textit{2D(b)} with regular and small diffusion.  
In Table \ref{cc-fa}, we present a summary for factor analysis of the two-species competition model. 
From the correlation matrix, we observed a strong correlation between final population density and diffusion rate and between final population density and growth and competition efficiency.

%\FloatBarrier
\begin{table}[h!]\centering
\scriptsize
\begin{tabular}{lrrrrrrr}\toprule
\textbf{Dimension} &\textbf{1D} &\textbf{2D(a)} &\textbf{2D(b)} &\textbf{1D} &\textbf{2D(a)} &\textbf{2D(b)} \\\cmidrule{1-7}
\textbf{Diffusion} &\textbf{Regular} &\textbf{Regular} &\textbf{Regular} &\textbf{Small} &\textbf{Small} &\textbf{Small} \\\cmidrule{1-7}
\textbf{Cum. Var.} &\textbf{60.37\%} &\textbf{63.71\%} &\textbf{67.37\%} &\textbf{67.18\%} &\textbf{67.53\%} &\textbf{67.20\%} \\\midrule
\textbf{Factor 1} &\cellcolor[HTML]{deeaf6} Diffusion 2 &\cellcolor[HTML]{deeaf6} Diffusion 2 &\cellcolor[HTML]{deeaf6} Diffusion 1 & Grow Compete 2 &Grow Compete 1 &Grow Compete 2 \\
\textbf{Factor 2} &\cellcolor[HTML]{deeaf6} Diffusion 1 &\cellcolor[HTML]{deeaf6} Diffusion 1 &\cellcolor[HTML]{deeaf6} Diffusion 2 & Grow Compete 1 &Grow Compete 2 &Grow Compete 1 \\
\textbf{Factor 3} & Grow Compete 1 &Grow Compete 1 & Grow Compete 2 &\cellcolor[HTML]{deeaf6} Diffusion 1 &\cellcolor[HTML]{deeaf6}Diffusion 2 &\cellcolor[HTML]{deeaf6}Diffusion 1 \\
\textbf{Factor 4} &Grow Compete 2 & Grow Compete 2 & Grow Compete 1 &\cellcolor[HTML]{deeaf6} Diffusion 2 &\cellcolor[HTML]{deeaf6}Diffusion 1 &\cellcolor[HTML]{deeaf6}Diffusion 2 \\
\bottomrule
\end{tabular}
\caption{Factor Analysis of two species system }
\label{cc-fa}
\end{table}
%\FloatBarrier

From the loading matrix of factor analysis, we observed that in a two-species system, when the diffusion rates of both species are regular, the dominant factor that causes variation is the diffusion rates of both species. The second most crucial factor is both species' growth and competition efficiency.
However, when the diffusion rates of both species are low, the dominant factor that causes variation becomes the growth and competition efficiency of both species instead of the diffusion rate. 
A possible explanation for the change of the dominant factor when the diffusion rate changes is the existence of the boundary conditions (set at zero). The effect of boundary is more severe when the diffusion rate is larger.

%\FloatBarrier
\begin{figure}[h!]
\centering
\vspace{10pt}
\begin{subfigure}{0.32\textwidth}
\centering
1D\\
\includegraphics[width=0.99\linewidth]{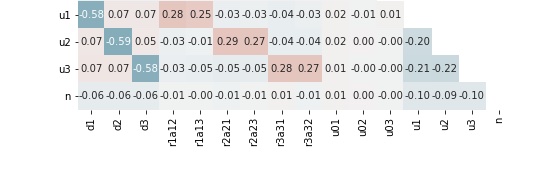}\\
\includegraphics[width=0.99\linewidth]{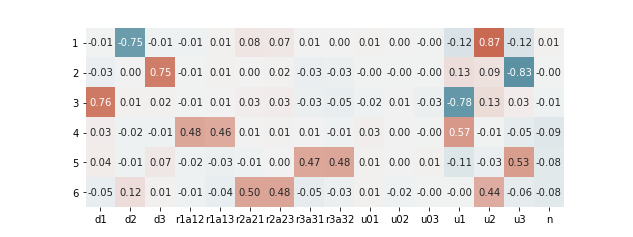}
\end{subfigure}
\begin{subfigure}{0.32\textwidth}
\centering
2D (a)\\
\includegraphics[width=1\linewidth]{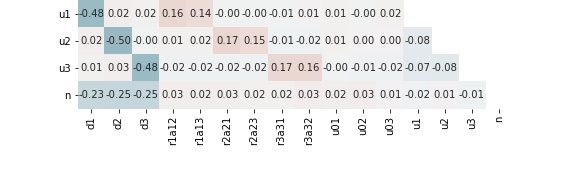}\\
\includegraphics[width=1\linewidth]{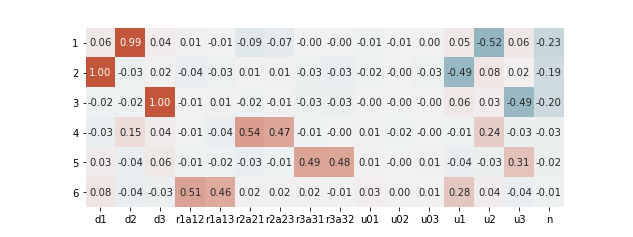}
\end{subfigure}
\begin{subfigure}{0.32\textwidth}
\centering
2D (b)\\
\includegraphics[width=1\linewidth]{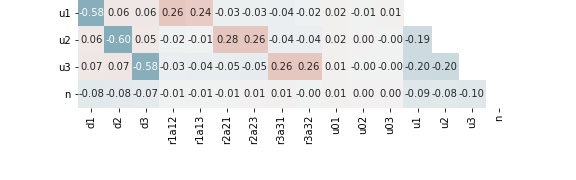}\\
\includegraphics[width=1\linewidth]{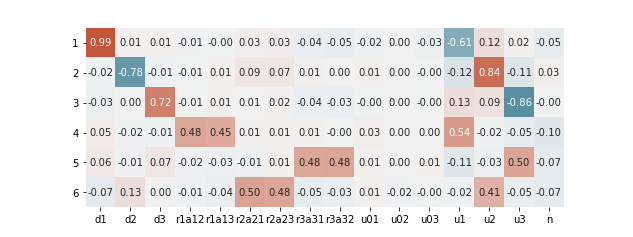}
\end{subfigure}\\
Regular diffusion
\\
\vspace{10pt}
\begin{subfigure}{0.32\textwidth}
\centering
\includegraphics[width=1\linewidth]{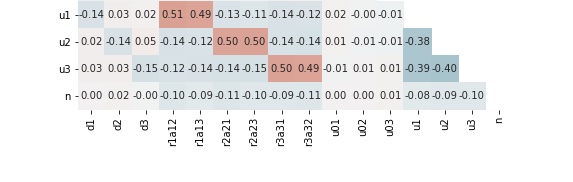}\\
\includegraphics[width=1\linewidth]{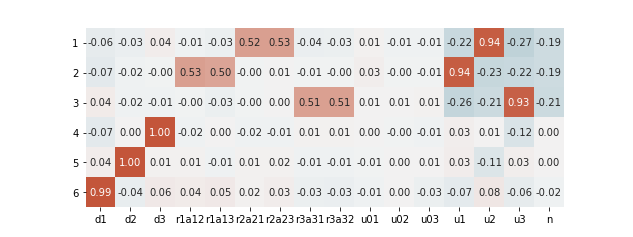}
\end{subfigure}
\begin{subfigure}{0.32\textwidth}
\centering
\includegraphics[width=1\linewidth]{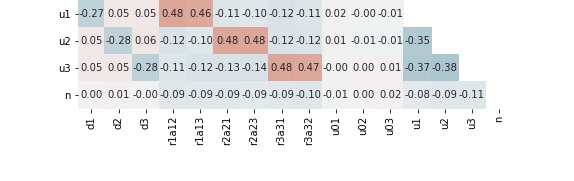}\\
\includegraphics[width=1\linewidth]{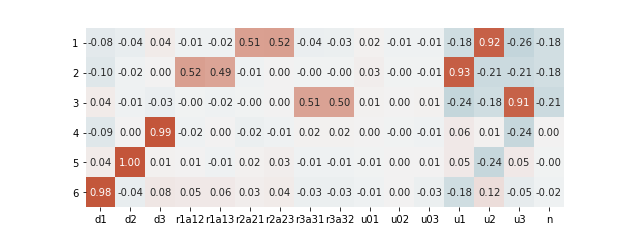}
\end{subfigure}
\begin{subfigure}{0.32\textwidth}
\centering
\includegraphics[width=1\linewidth]{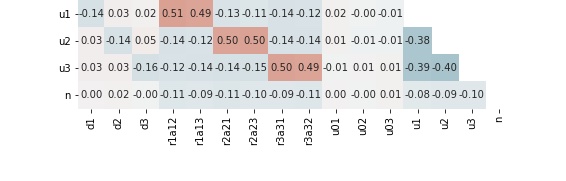}\\
\includegraphics[width=1\linewidth]{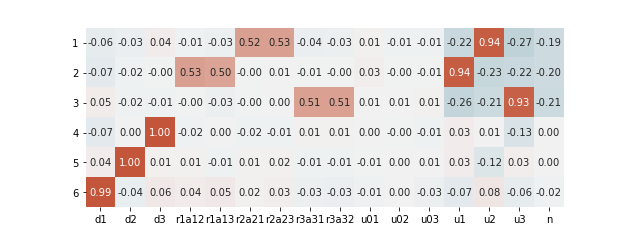}
\end{subfigure}
\\
Small diffusion
\caption{Three-species competition model. Correlation matrix (first row) and loading  (second row) }
\label{ccc-loading}
\end{figure}
%\FloatBarrier

%\FloatBarrier
\begin{table}[h!]\centering
\scriptsize
\begin{tabular}{lrrrrrrr}\toprule
\textbf{Dimension} &\textbf{1D} &\textbf{2D(a)} &\textbf{2D(b)} &\textbf{1D} &\textbf{2D(a)} &\textbf{2D(b)} \\\cmidrule{1-7}
\textbf{Diffusion} &\textbf{Regular} &\textbf{Regular} &\textbf{Regular} &\textbf{Small} &\textbf{Small} &\textbf{Small} \\\cmidrule{1-7}
\textbf{Cum. Var.} &\textbf{65.66\%} &\textbf{62.31\%} &\textbf{65.52\%} &\textbf{65.88\%} &\textbf{65.96\%} &\textbf{65.88\%} \\\midrule
\textbf{Factor 1} &\cellcolor[HTML]{deeaf6}Diffusion 2 &\cellcolor[HTML]{deeaf6}Diffusion 3 &\cellcolor[HTML]{deeaf6}Diffusion 1 &Grow Compete 2 &Grow Compete 2 &Grow Compete 2 \\
\textbf{Factor 2} &\cellcolor[HTML]{deeaf6}Diffusion 1 &\cellcolor[HTML]{deeaf6}Diffusion 2 &\cellcolor[HTML]{deeaf6}Diffusion 2 &Grow Compete 3 &Grow Compete 3 &Grow Compete 3 \\
\textbf{Factor 3} &\cellcolor[HTML]{deeaf6}Diffusion 3 &\cellcolor[HTML]{deeaf6}Diffusion 1 &\cellcolor[HTML]{deeaf6}Diffusion 3 &Grow Compete 1 &Grow Compete 1 &Grow Compete 1 \\
\textbf{Factor 4} &Grow Compete 2 &Grow Compete 3 &Grow Compete 2 &\cellcolor[HTML]{deeaf6}Diffusion 3 &\cellcolor[HTML]{deeaf6}Diffusion 2 &\cellcolor[HTML]{deeaf6}Diffusion 3 \\
\textbf{Factor 5} &Grow Compete 3 &Fin. Pop. 3 &Grow Compete 3 &\cellcolor[HTML]{deeaf6}Diffusion 2 &\cellcolor[HTML]{deeaf6}Diffusion 3 &\cellcolor[HTML]{deeaf6}Diffusion 2 \\
\textbf{Factor 6} &Grow Compete 1 &Grow Compete 2 &Grow Compete 1 &Grow Compete 1 &\cellcolor[HTML]{deeaf6}Diffusion 1 &Grow Compete 1 \\
\textbf{Factor 7} &Init. Pop. 3 &Grow Compete 1 &Init. Pop. 3 &\cellcolor[HTML]{deeaf6}Diffusion 1 &Grow Compete 1 &\cellcolor[HTML]{deeaf6}Diffusion 1 \\
\bottomrule
\end{tabular}
\caption{Factor Analysis of three species system}
\label{ccc-fa}
\end{table}
%\FloatBarrier

In Figure \ref{ccc-loading}, we present three - species competition system correlation matrices as well as loading matrix for each factor analysis for all three cases \textit{1D}, \textit{2D(a)} and \textit{2D(b)} with regular and small diffusion.  
In Table \ref{cc-fa}, we present a summary for factor analysis of the three-species competition model.
We observed that the three-species competing system is similar to the two-species competing system. There is a strong correlation between final population density and diffusion rate and between final population density and growth and competition efficiency. When the diffusion rate is regular, the dominant factors that cause variation are diffusion rates; in contrast, when the diffusion rates are low, the dominant factors become the growth and competition efficiency of the three species.

%\FloatBarrier
\begin{table}[h!]
\centering
\begin{tabular}{|c|cccc|}
 \hline
      & 00     & 01     & 10     & 11     \\ 
\hline
\multicolumn{5}{|c|}{Regular diffusion} \\   
\hline
1d    & 2796   & 3059   & 3064   & 1081   \\
      & 27.96\% & 30.59\% & 30.64\% & 10.81\% \\ \hline
2d(a) & 6888   & 1487   & 1489   & 136    \\
      & 68.88\% & 14.87\% & 14.89\% & 1.36\%  \\ \hline
2d(b) & 3125   & 2904   & 3043   & 928    \\
      & 31.25\% & 29.04\% & 30.43\% & 9.28\%  \\
\hline
\multicolumn{5}{|c|}{Small diffusion} \\    
\hline
1d    & 0      & 2251   & 2279   & 5470   \\
      & 0.00\%  & 22.51\% & 22.79\% & 54.70\% \\ \hline
2d(a) & 26     & 2522   & 2533   & 4919   \\
      & 0.26\%  & 25.22\% & 25.33\% & 49.19\% \\ \hline
2d(b) & 0      & 2184   & 2192   & 5624   \\
      & 0.00\%  & 21.84\% & 21.92\% & 56.24\%\\
\hline      
\end{tabular}
\caption{Survival Status for two species system}
\label{cc-survive}
\end{table}
%\FloatBarrier

%\FloatBarrier
\begin{table}[h!]
\centering
\begin{tabular}{|c|cccccccc|}
\hline
Group & 000    & 001    & 010    & 100    & 011    & 110    & 101    & 111    \\
\hline
\multicolumn{9}{|c|}{Regular diffusion} \\    
\hline
1d    & 1442   & 2026   & 2031   & 1987   & 806    & 822    & 730    & 156    \\
      & 14.42\% & 20.26\% & 20.31\% & 19.87\% & 8.06\%  & 8.22\%  & 7.30\%  & 1.56\%  \\ \hline
2d(a) & 5800   & 1262   & 1319   & 1234   & 131    & 130    & 118    & 6      \\
      & 58.00\% & 12.62\% & 13.19\% & 12.34\% & 1.31\%  & 1.30\%  & 1.18\%  & 0.06\%  \\ \hline
2d(b) & 1754   & 2012   & 2034   & 1985   & 712    & 728    & 651    & 124    \\
      & 17.54\% & 20.12\% & 20.34\% & 19.85\% & 7.12\%  & 7.28\%  & 6.51\%  & 1.24\%  \\
\hline
\multicolumn{9}{|c|}{Small diffusion} \\    
\hline
1d    & 0      & 639    & 648    & 556    & 2031   & 2023   & 1943   & 2160   \\
      & 0.00\%  & 6.39\%  & 6.48\%  & 5.56\%  & 20.31\% & 20.23\% & 19.43\% & 21.60\% \\ \hline
2d(a) & 3      & 858    & 835    & 773    & 1957   & 2001   & 1863   & 1710   \\
      & 0.03\%  & 8.58\%  & 8.35\%  & 7.73\%  & 19.57\% & 20.01\% & 18.63\% & 17.10\% \\ \hline
2d(b) & 0      & 651    & 655    & 566    & 2030   & 2023   & 1941   & 2134   \\
      & 0.00\%  & 6.51\%  & 6.55\%  & 5.66\%  & 20.30\% & 20.23\% & 19.41\% & 21.34\%\\
 \hline
\end{tabular}
\caption{Survival Status for three species system}
\label{ccc-survive}
\end{table}
%\FloatBarrier

For two-species competing system, we released control for all variables including diffusion rates, growth and competition efficiency, and initial population density, ran 10,000 simulations for all \textit{1D}, \textit{2D(a)}, and \textit{2D(b)} spatial boundary conditions, and summarized the proportion of different survival status groups, as is in Table \ref{cc-survive}. We observed that \textit{1D} and \textit{2D(b)} spatial boundary conditions give similar proportion of survival groups, while \textit{2D(a)} gives more different proportion. 
When there is regular diffusion, for \textit{1D} and \textit{2D(b)} spatial boundary conditions, there is around 30 percent of simulations that result in either 00(no species survives), 01 (species two survives), or 10 (species one survives). In comparison, only around 10 percent of simulations result in both species surviving together. For \textit{2D(a)} spatial boundary conditions, 68.88 percent of simulations result in 00 (both species go extinct), while around 14.87 percent and 14.89 percent result in 01 (species two survives) or 10 (species one survives), and very rare, about 1.36 percent of the simulation result in 11 (both species survive). To conclude, when there is regular diffusion, it's generally harder to arrive at a final solution where both species survive; it's even harder, almost impossible, for both species to survive when we are approximating a pond (\textit{2D(a)} scenario).

When there is small diffusion, all three spatial boundary conditions (\textit{1D}, \textit{2D(a)}, and \textit{2D(b)}) give a similar proportion of survival status. It is most likely that both species survive together (around 50 percent of simulations). However, it is almost impossible that both species will go extinct. Furthermore, it is equally likely that only one species survives (proportion of simulation both around 20-25 percent).

For three-species competing system, we released control for all variables including diffusion rates, growth and competition efficiency, and initial population density, ran 10,000 simulations for all \textit{1D}, \textit{2D(a)}, and \textit{2D(b)} spatial boundary conditions, and summarized the proportion of different survival status groups, as is in Table \ref{ccc-survive}. We observed that \textit{1D} and \textit{2D(b)} spatial boundary conditions give similar proportion of survival groups, while \textit{2D(a)} gives more different proportion.  
When there is regular diffusion, the probability of one species survives (001, 010, and 100) is similar to one another (around 10 percent for \textit{1D}/\textit{2D(b)} scenario, and around 12 percent for \textit{2D(a)} scenario), while the probability of two species survive (011, 110 and 101) is similar to one another(around 12 percent for \textit{1D}/\textit{2D(b)} scenario, and around 1 percent for \textit{2D(a)} scenario). In general, when there are three species, it is almost impossible for all three species to survive together and no species go extinct; what's worse, for \textit{2D(a)} scenario which approximate the lake, around 58 percent of simulation result in 000 (no species survive). 
When there is small diffusion, all three spatial boundary conditions (\textit{1D}, \textit{2D(a)}, and \textit{2D(b)}) give a similar proportion of survival status. It's equally likely that two species survive or all species survive (the proportion for survival group 011/110/101/111 are all around 20 percent), and it's equally likely that only one species survive (the proportion for survival group 001/010/100 are all around 6-8 percent), and it's almost impossible that all three species go extinct.

\section{Conclusion}

A mathematical model of the spatial-temporal multi-species competition is considered.  
A discrete system is constructed using a finite volume approximation with a semi-implicit time approximation.  
The numerical results for two- and three-species models are presented for several exceptional cases of the parameters related to the survival status.  We considered the one and two-dimensional model problems with two cases of the boundary conditions for the two-dimensional case.  
First, the effect of diffusion is investigated numerically. In special cases where parameters are fixed, we observed that the effect of boundary constraint is more severe in regular diffusion groups than in small diffusion groups, causing a lower population density both in the middle and near the boundary of the domain. %We also observed that the dynamic of \textit{1D} is similar to \textit{2D(b)}, while \textit{2D(a)} gives different dynamic. This suggest that we can use \textit{1D} to approximate \textit{2D(b)}, so as to save computation time. What is more, 
Furthermore from general cases where we release holds of diffusion while keeping other parameters fixed, we observed that different combination of diffusion rates of species in the system leads to different final survival status of species. When the combination of both species' diffusion rates is on the border of different survival groups, it takes more time for these groups of species to reach equilibrium.
Second, the effect of the initial condition is investigated numerically. We release holds of initial conditions while keeping other parameters fixed. We observed similar patterns for both two-species competing systems and three-species competing for systems. Take a two-species competing system as an example; we observed that disregarding the value of the initial population of both species, the final equilibrium population will rest at a point of focus. The farther the combination point of initial population density is from the final equilibrium combination point, the more time steps it would take to reach the equilibrium.

Finally, the impact of parameters on the system stability is considered by simulating the spatial-temporal model with random input parameters. Factor analysis and statistical summary of the survival status of species in the system was performed.
%At last, we release the holds for all parameters and simulate 10,000 results for random parameters under regular and small diffusion rates, then perform factor analysis and statistical summary of the survival status of species in the system. In the factor analysis, w
We observed %a similar behavior of the results for two and three-species competition models. D
that diffusion rates are the dominant factor. when diffusion rates are regular and on the same scale as other parameters. In contrast, when diffusion rates are small, which are ten times smaller in the scale of other parameters, growth and competition rates become the dominant factors. 
In a statistical summary of species survival status, we observed a similar pattern for both two-species competing systems and three-species competing systems. %Take two-species competing system as an example, when there is small diffusion, all three spatial boundary conditions (\textit{1D}, \textit{2D(a)}, and \textit{2D(b)}) give a similar proportion of survival status. It is most likely that both species survive together. However, it is almost impossible that both species will go extinct. Furthermore, it is equally likely that only one species survives.
%However, in regular diffusion, we observed that \textit{1D} and \textit{2D(b)} spatial boundary conditions give similar proportion of survival groups, while \textit{2D(a)} gives more different proportion of species survival status. Take a two-species competing system as an example; when there is regular diffusion, it is generally harder to arrive at a final solution where both species survive; it is even harder, almost impossible, for both species to survive when we are approximating a pond (\textit{2D(a)} scenario).

\bibliographystyle{plain}
\bibliography{multispecies_cite}

\end{document}